\documentclass{amsart}

\usepackage{amssymb,amsmath,amsthm}
\usepackage{enumerate}
\usepackage{url}
\usepackage{tikz-cd}
\usepackage{mathtools}
\usepackage{hyperref}

\DeclareMathOperator{\G}{G}
\DeclareMathOperator{\Gal}{Gal}
\DeclareMathOperator{\GL}{GL}

\DeclareMathOperator{\ord}{ord}

\DeclareMathOperator{\Aut}{Aut}
\DeclareMathOperator{\residue}{Res}
\DeclareMathOperator{\Rad}{Rad}
\DeclareMathOperator{\degree}{deg}

\theoremstyle{plain}
\newtheorem{theorem}{Theorem}[section]

\newtheorem{corollary}[theorem]{Corollary}
\newtheorem{proposition}[theorem]{Proposition}

\theoremstyle{definition}

\newtheorem{example}[theorem]{Example}
\newtheorem{remark}[theorem]{Remark}

\newcommand{\Z}{\mathbb{Z}}
\newcommand{\Q}{\mathbb{Q}}

\newcommand{\F}{\mathbb{F}}
\newcommand{\Proj}{\mathbb{P}}

\title{Explicitly bounding perfect powers in elliptic divisibility sequences}
\author{Sander R. Dahmen}
\author{Joey M. van Langen}
\address{
Department of Mathematics,
Vrije Universiteit Amsterdam,
De Boelelaan 1111,
1081 HV Amsterdam,
The Netherlands}
\email{s.r.dahmen@vu.nl}
\email{j.m.van.langen@outlook.com}
\subjclass[2010]{11G05, 11D41}
\keywords{Modular form, Diophantine equation, Elliptic divisibility sequence}
\date{\today}
\thanks{Supported by NWO Vidi grant 639.032.613}

\begin{document}

\begin{abstract}
  In this paper we consider elliptic divisibility sequences generated by a point on an elliptic curve over $\Q$ with $j$-invariant $1728$ given by an integral short Weierstrass equation.
  For several different such elliptic divisibility sequences, we determine explicitly a finite set of primes such that for all primes $l$ outside this set, the elliptic divisibility sequence contains no $l$-th powers.
  Our approach uses, amongst other ingredients, the modular method for $\Q$-curves.
  The corresponding explicit computations fit into a general computational framework for $\Q$-curves, for which they provide illustrative examples.
\end{abstract}

\maketitle

\section{Introduction}

Given a non-singular Weierstrass equation 
\begin{equation} \label{we}
E: y^2+a_1xy+a_3y=x^3+a_2x^2+a_4x+a_6
\end{equation}
with integral coefficients and a non-torsion point $P \in E(\mathbb{Q})$, 
one can use the group law on $E$ to write for a positive integer $m$
\begin{equation} \label{B_m}
x(mP) = \frac{A_m}{B_m^2}
\end{equation}
for coprime integers $A_m, B_m$ with $B_m>0$.
The sequence $\left(B_m\right)_{m \in \Z_{>0}}$ is known as an elliptic divisibility sequence. In particular, as is well known, for all natural numbers $n, m$ it satisfies
\begin{equation}\label{eqn:edsp}
n|m \Rightarrow B_n|B_m.
\end{equation}
Some
finiteness results have been achieved for perfect powers in elliptic
divisibility sequences; see e.g.~\cite{EverestReynoldsStevens07} and~\cite{Reynolds12}. In particular, if $E$ is a Mordell curve (so with $j$-invariant $j(E)=0$) and
$B_1>1$, then there are finitely many perfect powers in $(B_m)$; see \cite[Theorem 1.2]{Reynolds12}.
The method of proof includes a (regular) modular approach using elliptic curves over $\Q$.

Similarly, if $E$ has $j$-invariant $1728$ and $B_1>1$, then
Dahmen and Reynolds study (in a preprint) again the finiteness of perfect powers in the corresponding sequence.
In this case, a Frey $\Q$-curve is associated to the problem.
Its field of definition can sometimes taken to be $\Q$, but \lq generically\rq\ it is a quadratic number field. 
We do not strictly depend on their results, but our Frey $\Q$-curve construction is similar.
If in the Weierstrass equation~\ref{we}, we have that $a_4$ is a positive integer and all other $a_i$ are zero,
then the corresponding Frey $\Q$-curve is defined over a real quadratic number field (or sometimes already over $\Q$),
hence one can also try to perform a Hilbert modular approach in this situation.
This has in fact been done in~\cite{Alfaraj23}, showing the finiteness of perfect powers if $B_1>1$ (still assuming $a_4>0$).

In the remainder of this paper, we focus on the $j=1728$ case. The aim is to show that for many choices of $E,P$ we can explicitly find a finite set $S$ of primes such for all primes $l$ not in $S$ the associated elliptic divisibility sequence contains no $l$-th powers. We will work out several examples, illustrating the power of the method.
Our main Diophantine results are given in Section~\ref{sec:EDS examples}, namely Theorem~\ref{thm:D125}, Theorem~\ref{thm:Dm17}, as well as the results mentioned in Subsection~\ref{FurtherExamples}. Most notably are cases (ii), (iii), and (viii) in Table~\ref{tab:EDSresults}, as no alternative approaches for these seem to be available in the literature. In principle, we believe that many more examples could be computed along similar lines of our approach.

In order to compute the levels of newforms associated with a Frey
$\Q$-curve~${ E }$, one needs to compute a particular twist of the
curve. Given a Galois number field~${ K }$ over which~${ E }$ is
completely defined and for which a splitting map factors
over~${ G_\Q^K }$ (see Section~\ref{sec:modularity} for definitions),
this twist can be computed from a map~${ \alpha : G_\Q^K \to K^* }$
whose coboundary is a~${ 2
}$-cocycle~${ \left(G_\Q^K\right)^2 \to \{ \pm 1 \} }$ associated with
the Frey~${ \Q }$-curve. In the literature, a
map~${ \alpha : G_\Q^K \to \mathcal{O}_K^* }$ often suffices, but we
shall show that this is not the case in the example corresponding to
Theorem~\ref{thm:Dm17}. We will demonstrate that we can instead
take~${ \alpha : G_\Q^K \to \mathcal{O}_{K, S}^* }$
where~${ \mathcal{O}_{K, S}^* }$ is the ring of~${ S }$-units for a
finite set~${ S }$.

The underlying computations in this paper have been automated by the second-named
author as part of the framework~\cite{framework} for computing with (Frey) $\Q$-curves.
For more details on this general framework, see \cite{LangenThesis}.
The specific computations for this paper can be found at:
\begin{center}
\url{https://github.com/sanderdahmen/Bounding-perfect-powers-in-EDS},
\end{center}
where we note (again) that these are based on the framework~\cite{framework}.

\section{Associating a $\mathbb{Q}$-curve}

Let $D$ be a nonzero integer and consider the Weierstrass equation
\begin{equation}\label{eqn j=1728 curve}
E_D: y^2=x^3+Dx.
\end{equation}
Note that~${ T = (0, 0) \in E_D(\Q) }$ is a~${ 2 }$-torsion point.
For any point~${ P \in E_D(\Q) \setminus \{ \mathcal{O}, T \} }$
take~${ \hat{P} = T - P \in E_D(\Q) \setminus \{ \mathcal{O}, T \}
}$. Since~${ P + \hat{P} + T = \mathcal{O} }$, all of these points are
on the line~${ y = \frac{y_P}{x_P} x }$. Substituting this in the
equation of~${ E_D }$ and comparing coefficients of $x^2$, tells us that
\begin{equation*}
  x_P + x_{\hat{P}} + 0 = \left( \frac{y_P}{x_P} \right)^2 = \frac{x_P^3 + D x_P}{x_P^2} = x_P + \frac{D}{x_P},
\end{equation*}
hence~${ x_P x_{\hat{P}} = D }$. If we write in lowest terms,
\begin{equation*}
  P = \left( \frac{A}{B^2}, \frac{C}{B^3} \right) \text{ and }
  \hat{P} = \left( \frac{\hat{A}}{\hat{B}^2}, \frac{\hat{C}}{\hat{B}^3} \right),
\end{equation*}
so $\gcd(AC,B)=\gcd(\hat{A}\hat{C},\hat{B})=1$ and $B,\hat{B}>0$,
then this implies~${ A \hat{A} = D B^2 \hat{B}^2 }$.
Since~${ \gcd(A, B) = \gcd(\hat{A}, \hat{B}) = 1 }$, there
are~${ a, \hat{a} \in \Z }$ such that~${ A = a \hat{B}^2 }$
and~${ \hat{A} = \hat{a} B^2 }$.  Note that~${ D = a \hat{a} }$ and
that~${ \gcd(A, D) = | a | \gcd(\hat{B}^2, \hat{a}) = | a | }$, so one
can compute~${ a }$ and~${ \hat{a} }$ easily, noting that the sign
of~${ a }$ must be the same as the sign of~${ A }$. Now substituting
the coordinates of~${ P }$ into the equation of~${ E_D }$ and
multiplying by~${ B^6 }$, tells us that
\begin{equation*}
  C^2 = a^3 \hat{B}^6 + a^2 \hat{a} B^4 \hat{B}^2.
\end{equation*}
Therefore, ${ C = a \hat{B} w }$ for some~${ w \in \Z }$, and dividing the equation above by~${ a^2 \hat{B}^2 }$ gives
\begin{equation}
  \label{eqn basic 3-term}
  w^2 = a \hat{B}^4 + \hat{a} B^4.
\end{equation}
Note that in this equation~${ B }$ and~${ \hat{B} }$ are coprime
as~${ \gcd(B, \hat{B}) \mid \gcd(B, A) = 1 }$. Furthermore, substituting the
coordinates of~${ \hat{P} }$ into the equation for~${ E_D }$ gives rise
to the same Diophantine equation
as~${ \hat{C} = \hat{B}^3 \frac{y_P}{x_P} x_{\hat{P}} = \hat{a} B w
}$.
For future reference, we note that actually $B,\hat{B}$, and $w$ are pairwise coprime since also $\gcd(w,B)|\gcd(C,B)=1$ and $\gcd(w,\hat{B})|\gcd(\hat{C},\hat{B})=1$.

The condition that $B$ is a perfect (say $l$-th) power now leads to a
so-called generalized Fermat equation of signature $(2,4,4l)$, hence also of signature $(2,4,l)$. To the
Fermat equation $y^2 - d x^4 = e z^l$ (for given nonzero integers
$d, e$), we associate (basically as in \cite{DieulefaitUrroz09}) the
Frey $\Q$-curve
\begin{equation}\label{eqn E a,b,d,gamma}
E_{d,x,y}: Y^2=X^3+4 \sqrt{d} x X^2+2(d x^2 + \sqrt{d} y) X.
\end{equation}
We associate to~\eqref{eqn basic 3-term} the $\Q$-curve above with
$d=a$, $x = \hat{B} =: z $, and $y=w$, and for later reference a \lq
twisting parameter\rq\ $\gamma$ (some algebraic integer), i.e.
\begin{equation}
  \label{eqn E gamma}
  E^{\gamma}_{a, z, w} : Y^2=X^3+ 4 \sqrt{a} z \gamma X^2+2\left(a z^2 + \sqrt{a} w \right) \gamma^2 X.
\end{equation}
This model has~${c_4}$-invariant
\begin{equation*}
  c_{4, a, z, w}^\gamma = - 2^5 \sqrt{a} \left( 3 w - 5 \sqrt{a} z^2 \right) \gamma^2
\end{equation*}
discriminant
\begin{equation*}
\Delta_{a, z, w}^{\gamma} = - 2^9 \sqrt{a}^3 \left( w - \sqrt{a} z^2 \right)
\left( w + \sqrt{a} z^2 \right)^2 \gamma^6,
\end{equation*}
and~${ j }$-invariant
\begin{equation*}
  j_{a, z, w} = 2^6 \frac{(3 \, w - 5 \, \sqrt{a} z^2)^3}{(w + \sqrt{a} z^2)^2 (w - \sqrt{a} z^2)}.
\end{equation*}

\begin{proposition}
  \label{thm:c4Deltacoprime}
  The invariants~${ c_{4, a, z, w}^\gamma }$
  and~${ \Delta_{a, z, w}^\gamma }$ are coprime outside all primes
  dividing~${ 2 \gamma a }$.
\end{proposition}
\begin{proof}
  Suppose a finite prime~${ \mathfrak{p} }$
  divides~${ c_{4, a, z, w}^\gamma }$
  and~${ \Delta_{a, z, w}^\gamma }$ but not~${ 2 \gamma a }$,
  then~${ \mathfrak{p} }$ should divide~${ z }$ and~${ w }$. Note
  however that~${ w \mid \hat{C} }$ and~${ z = \hat{B} }$,
  so~${ \mathfrak{p} \mid \gcd(\hat{C}, \hat{B}) = 1 }$, a
  contradiction.
\end{proof}
\begin{corollary}
  \label{thm:nonadditive}
  At all finite primes not dividing~${ 2 \gamma a }$ the
  model~${ E_{a, z, w}^\gamma }$ is minimal and has semi-stable (i.e. non-additive)
  reduction.
\end{corollary}

\begin{proposition}
  \label{thm:potentialmultiplicative}
  Suppose~${ B }$ is divisible by a prime number~${ p }$,
  with~${ p^3 \mid B }$ if~${ p = 2 }$, then~${ j_{a, z, w} }$ is not
  integral at primes above~${ p }$. In
  particular, ${ E_{a, z, w}^{\gamma} }$ has potentially multiplicative
  reduction at such primes.
\end{proposition}
\begin{proof}
  Suppose that~${ p \mid B }$, then we know that~${ p \nmid A = a z^2 }$
  and~${ p \nmid C = a z w }$.
  Since ${ \hat{a} B^4 = (w + \sqrt{a} z^2) (w - \sqrt{a} z^2) }$, we get that primes
  above~${ p }$ divide the denominator, but not the numerator
  of~${ j_{a, z, w} }$ if~${ p \ne 2 }$.

  When~${ p = 2 }$, let~${ \mathfrak{p} }$ be a prime above~${ 2 }$
  and suppose~${ \ord_{\mathfrak{p}} j_{a, z, w} \ge 0 }$. Writing
  this out gives
  \begin{equation}
    \label{eqn j integral}
    6 \ord_{\mathfrak{p}} 2 + 3 \ord_{\mathfrak{p}}(3 \, w - 5 \sqrt{a} z^2) \ge 2 \ord_{\mathfrak{p}}(w + \sqrt{a} z^2) + \ord{\mathfrak{p}}(w - \sqrt{a} z^2).
  \end{equation}
  If~${ \ord_{\mathfrak{p}}(w - \sqrt{a} z^2) >
    \ord_{\mathfrak{p}} 2 }$ then~${ \ord_{\mathfrak{p}}(3 \, w - 5 \sqrt{a} z^2) =
    \ord_{\mathfrak{p}}(w + \sqrt{a} z^2) = \ord_{\mathfrak{p}} 2
  }$ and
  hence~${ \ord_{\mathfrak{p}}(w - \sqrt{a} z^2) =
    \ord_{\mathfrak{p}} 2 \left( \ord_2 \hat{a} + 4 \ord_2 B - 1
    \right) }$. Substituting this in equation~\eqref{eqn j integral}
  gives us
  \begin{equation*}
    9 \ge \ord_2 \hat{a} + 4 \ord_2 B + 1,
  \end{equation*}
  contradicting~${ 8 \mid B }$.  We thus
  have~${ \ord_{\mathfrak{p}}(w - \sqrt{a} z^2) = \ord_{\mathfrak{p}}
    2 }$ and therefore we
  have that~${ \ord_{\mathfrak{p}}(w + \sqrt{a} z^2) = \ord_{\mathfrak{p}}
    2 \left( \ord_2 \hat{a} + 4 \ord_2 B - 1 \right) > 3
    \ord_{\mathfrak{p}} 2 }$. We thus find that
  \begin{equation*}
    3 \ord_{\mathfrak{p}}(3 \, w - 5 \sqrt{a} z^2) \ge \ord_{\mathfrak{p}} 2 \left( 2 \ord_2 \hat{a} + 8 \ord_2 B - 1 \right) > 9 \ord_{\mathfrak{p}} 2.
  \end{equation*}
  This contradicts that~${ 3 (w + \sqrt{a} z^2) - (3 \, w - 5 \sqrt{a} z^2) = 8 w }$ has
  valuation~${ 3 \ord_{\mathfrak{p}} 2 }$, thereby completing the proof.
\end{proof}
\begin{corollary}
  \label{thm:nonCM}
  If~${ B }$ is divisible by a prime number~${ p }$,
  with~${ p^3 \mid B }$ if~${ p = 2 }$, then~${ E_{a, z, w} }$ does
  not have complex multiplication.
\end{corollary}
\begin{proof}
  Follows directly from~\cite[II, Theorem~6.1]{Silverman94}
\end{proof}

For some of the $\Q$-curve computations~\cite{PacettiVillagra22, PacettiVillagra23} could also apply.

Whenever the point~${ P }$ is not clear from the context, we will add
an index~${ P }$ to the variables mentioned above,
e.g.~${ a_P, A_P, }$ etc. If $E_D$ has positive rank and
$P \in E_D(\Q)$ is a non-torsion point we will denote the variables
associated to~${ P_m := m P }$ by an index~${ m \in \Z_{>0} }$,
e.g.~${ a_m := a_{P_m}, A_{m} := A_{P_m}}$, etc.

\subsection{$a$ only depends on class modulo $[2] E_D(\Q)$}
\label{sec:a depends on 2E}
We will now show that the value of~${ a }$ for a point only depends on
its class modulo the image of the multiplication by~2 map. This will
in particular show that the Frey~${ \Q }$-curve associated
with~${ P_m }$ (${ m \in \Z_{>0}}$, ${ P \in E_D(\Q) }$ non-torsion)
only depends on the parity of~${ m }$. Therefore, the construction~\eqref{eqn E gamma} yields
at most~${ 2 }$ distinct Frey~${ \Q }$-curves up to twisting by $\gamma$ (in the unknowns $z$ and $w$) associated to an elliptic
divisibility sequence.

Note that the short exact sequence
\begin{equation*}
 \begin{tikzcd}
    1 \arrow[r] & E_D[2] \arrow[r] & E_D(\overline{\Q}) \arrow[r, "{[2]}"] & E_D(\overline{\Q}) \arrow[r] & 1,
  \end{tikzcd}
\end{equation*}
induces a long exact sequence in Galois cohomology containing a map
\begin{equation*}
  \alpha : E_D(\Q) = H^0 \left( \Gal\left(\overline{\Q} / \Q \right), E_D(\overline{\Q}) \right) \to H^1 \left( \Gal\left(\overline{\Q} / \Q \right), E_D[2] \right).
\end{equation*}
The kernel of this map is precisely the image
of~${ [2]: E_D(\Q) \to E_D(\Q) }$, hence this map will be useful for proving
properties of points modulo this image. To work with the codomain of
this map, we note that~${ E_D[2] = \{ \mathcal{O}, T, T_{+1}, T_{-1} \} }$, where
\begin{equation*}
  T_{\delta} = (\delta \sqrt{-D}, 0),
\end{equation*}
for a fixed choice of~${ \sqrt{-D} }$.

For a fixed~${ [f] \in H^1(G_\Q, E_D[2]) }$ we can choose
functions~${ f_\delta : G_\Q \to \Z / 2\Z }$ such that
\begin{equation*}
  f(\sigma) = f_{+1}(\sigma) T_{+1} + f_{-1}(\sigma) T_{-1} \text{ for all } \sigma \in G_\Q.
\end{equation*}
Since~${ f }$ is a cocycle we find
that~${ f(\sigma) + \prescript{\sigma}{}f(\tau) = f(\sigma \tau) }$
for all~${ \sigma, \tau \in G_\Q }$, hence
\begin{equation*}
  \left\{
    \begin{array}{ll}
      f_{+1}(\sigma) + f_{+1}(\tau) = f_{+1}(\sigma \tau) \text{ and }
      f_{-1}(\sigma) + f_{-1}(\tau) = f_{-1}(\sigma \tau) &
                                                           \text{ if } \sigma \in G_{\Q(\sqrt{-D})} \\
      f_{+1}(\sigma) + f_{-1}(\tau) = f_{+1}(\sigma \tau) \text{ and }
      f_{-1}(\sigma) + f_{+1}(\tau) = f_{-1}(\sigma \tau) & \text{ otherwise.}
    \end{array}
  \right.
\end{equation*}
This implies
that~${ f_{\delta}|_{G_{\Q(\sqrt{-D})}} \in \hom(G_{\Q(\sqrt{-D})}, \Z /
  2\Z) \cong \Q\left(\sqrt{-D}\right)^* /
  \left(\Q\left(\sqrt{-D}\right)^*\right)^2 }$, hence we can construct
maps~${ \alpha_\delta : E_D(\Q) \to \Q\left(\sqrt{-D}\right)^* /
  \left(\Q\left(\sqrt{-D}\right)^*\right)^2 }$, by combining the
map~${ \alpha }$ with~${ [f] \mapsto f_\delta|G_{\Q(\sqrt{-D})}
}$. Note that the latter is well-defined as the relevant coboundaries
are
\begin{equation*}
  \left\{
    \begin{array}{ll}
      \partial \mathcal{O} = \partial T & = \left(\sigma \mapsto \mathcal{O}\right) \\
      \partial T_{+1} = \partial T_{-1} & = \left(\sigma \mapsto
                                      \left \{
                                      \begin{array}{ll}
                                        \mathcal{O} & \text{ if } \sigma \in G_{\Q(\sqrt{-D})} \\
                                        T & \text{ otherwise}
                                      \end{array}
                                      \right.
                                      \right).
    \end{array}
  \right.
\end{equation*}
In case~${ -D }$ is a square it is clear that the
maps~${ \alpha_\delta }$ encode all information of the
map~${ \alpha }$. In the case that there exists
a~${ \sigma \in G_{\Q} \setminus G_{\Q(\sqrt{-D})} }$ this remains
true and in fact one of the maps~${ \alpha_{\delta} }$ becomes
redundant. This can be seen as for
any~${ \tau \in G_{\Q(\sqrt{-D})} }$ we have
\begin{equation*}
  \begin{array}{l}
    f_{-1}(\tau) = f_{-1}(\tau \sigma^{-1}) - f_{-1}(\sigma^{-1})
    = f_{+1}(\sigma \tau \sigma^{-1}) - f_{+1}(\sigma) - f_{-1}(\sigma^{-1}) \\
    \phantom{f_{-1}(\tau)} = f_{+1}(\sigma \tau \sigma^{-1}) - f_{+1}(1)
    = f_{+1}(\sigma \tau \sigma^{-1}).
  \end{array}
\end{equation*}
Furthermore, the value of~${ f_{+1} }$ at any element
of~${ G_\Q \setminus G_{\Q(\sqrt{-D})} }$ can be derived
from~${ f_{+1}(\sigma) }$. Since there exists a non-trivial coboundary
in this case, the exact value of~${ f_{+1}(\sigma) }$ does not matter
for the cohomology class~${ [f] }$.

By making the maps~${ \alpha_\delta : E_D(\Q) \to \Q\left(\sqrt{-D}\right)^* /
  \left(\Q\left(\sqrt{-D}\right)^*\right)^2 }$ explicit, we obtain the following result.
\begin{proposition}
  \label{thm:descent}
  There exist
  homomorphisms
  \begin{equation*}
    { \alpha_{\delta} : E_D(\Q) \to
      \Q\left(\sqrt{-D}\right)^* /
      \left(\Q\left(\sqrt{-D}\right)^*\right)^2 }
  \end{equation*}  
  for~${ \delta \in \{ \pm 1 \} }$ such that the map
  \begin{equation*}
    (\alpha_{+1}, \alpha_{-1}) : E_D(\Q) \to \left( \Q\left(\sqrt{-D}\right)^* /
      \left(\Q\left(\sqrt{-D}\right)^*\right)^2\right)^2
  \end{equation*}
  has kernel~${ \{ [2] P : P \in E_D(\Q) \} }$.
  Explicitly, these homomorphisms are given by
  \begin{equation*}
    \alpha_\delta(P) = \left\{
      \begin{array}{ll}
        [x + \delta \sqrt{-D}] & \text{ if } P = (x, y) \\
        {[1]} & \text{ if } P = \mathcal{O}.
      \end{array}
    \right.
  \end{equation*}
\end{proposition}
\begin{proof}
  The existence of these maps has already been proven using the
  map~${ \alpha }$. To make them explicit, note that~${ \alpha }$ maps
  a point~${ P \in E_D(\Q) }$ to the class
  of~${ f : \sigma \mapsto \prescript{\sigma}{}Q - Q }$ for
  any~${ Q \in E_D(\overline{\Q}) }$ such that~${ [2] Q = P }$.

  Let~${ K = \Q(\sqrt{-D}) }$ and pick~${ \delta \in \{ \pm 1 \} }$
  and a point~${ P = (x, y) \in E_D(\Q) }$ arbitrary. Note
  that~${ \alpha_\delta(P) = [\gamma] }$ for some~${ \gamma \in K^* }$
  such that
  \begin{equation*}
    \left\{
      \begin{array}{ll}
        \alpha(\sigma) \in \{ \mathcal{O}, T_{-\delta} \} & \text{ if } \sigma \in G_{K(\sqrt{\gamma})} \\
        \alpha(\sigma) \in \{ T, T_{\delta} \} & \text{ if } \sigma \in G_K \setminus G_{K(\sqrt{\gamma})}.
      \end{array}
    \right.
  \end{equation*}
  All that remains to check is that in
  fact~${ \gamma = x + \delta \sqrt{-D} }$ satisfies this condition.

  Let~${ Q = (x_Q, y_Q) \in E_D(\overline{\Q}) }$ be a point such
  that~${ [2] Q = P }$. The tangent line to~${ E_D }$ at~${ Q }$ is
  given by
  \begin{equation*}
    Y = \frac{3 x_Q^2 + D}{2 y_Q} (X - x_Q) + y_Q.
  \end{equation*}
  By substituting this in the equation for~${ E }$, we get a cubic
  polynomial in~${ X }$ of which the roots are~${ x }$ and~${ x_Q }$
  twice. Looking at the coefficient of~${ X^2 }$ this tells us that
  \begin{equation*}
    2 x_Q + x = \left(\frac{3 x_Q^2 + D}{2 y_Q}\right)^2 = \frac{9 x_Q^4 + 6 D x_Q^2 + D^2}{4 x_Q^3 + 4 D x_Q },
  \end{equation*}
  hence~${ x_Q }$ is a root of the polynomial
  \begin{equation*}
    f(X) = X^4 - 4 \, x X^3 - 2 \, D X^2 - 4 \, x D X + D^2.
  \end{equation*}
  Note that all four roots of this polynomial correspond to the
  $x$-coordinates of the four possible points~${ Q }$, i.e.
  the~${ x }$-coordinates of~${ Q, Q + T, Q + T_{+1} }$,
  and~${ Q + T_{-1} }$.

  Now over the field~${ K(\sqrt{\gamma}) }$ the polynomial~${ f(X) }$
  splits as
  \begin{equation*}
    f(X) = (X^2 - \omega X - \delta \sqrt{-D} \omega - D)
    (X^2 - \overline{\omega} X - \delta \sqrt{-D} \overline{\omega} - D)
  \end{equation*}
  where~${ \omega = 2 x + 2 \frac{y}{\sqrt{\gamma}} }$
  and~${ \overline{\omega} = 2 x - \frac{y}{\sqrt{\gamma}}
  }$. Assuming without loss of generality that~${ x_Q }$ is a root of
  the first factor, the other root of the first factor
  is~${ \omega - x_Q }$, which we claim to be the~${ x }$-coordinate
  of~${ Q + T_{-\delta} }$. This would imply
  that 
  \begin{equation*}
    \left\{
      \begin{array}{ll}
        \prescript{\sigma}{}Q - Q \in \{ \mathcal{O}, T_{-\delta} \} & \text{ if } \sigma \in G_{K(\sqrt{\gamma})} \\
        \prescript{\sigma}{}Q - Q \in \{ T, T_{\delta} \} & \text{ if } \sigma \in G_K \setminus G_{K(\sqrt{\gamma})},
      \end{array}
    \right.
  \end{equation*}
  as only~${ \sigma \not\in G_{K(\sqrt{\gamma})} }$ switch the factors
  of~${ f(X) }$.

  It thus remains to prove the claim that~${ \omega - x_Q }$ is
  the~${ x }$-coordinate of~${ Q + T_{-\delta} }$. Note that the line
  through~${ Q }$ and~${ T_{-\delta} }$ is given by
  \begin{equation*}
    Y = \frac{y_Q}{x_Q + \delta \sqrt{-D}} (X + \delta \sqrt{-D}).
  \end{equation*}
  Doing the same substitution trick, we find that
  the~${ x }$-coordinate of~${ Q + T_{-\delta} }$ is
  \begin{equation*}
    \begin{array}{rl}
      \beta & = \left(\frac{y_Q}{x_Q + \delta \sqrt{-D}}\right)^2 - x_Q + \delta \sqrt{-D} \\
            & = \frac{x_Q (x_Q - \delta \sqrt{-D})}{x_Q + \delta \sqrt{-D}} - (x_Q - \delta \sqrt{-D}) \\
            & = \frac{-\delta \sqrt{-D} (x_Q - \delta \sqrt{-D})}{x_Q + \delta \sqrt{-D}}.
    \end{array}
  \end{equation*}
  To see that this is the same as~${ \omega - x_Q }$, note that
  \begin{equation*}
  \begin{array}{rl}
    (\beta + x_Q - 2 \, x)^2
    & = \left(\frac{x_Q^2 - D}{x_Q - \delta \sqrt{-D}} - 2 \, x \right)^2 \\
    & = \frac{x_Q^4 - 2 D x_Q^2 + D^2}{(x_Q - \delta \sqrt{-D})^2} - 4 \, x \frac{x_Q^2 - D}{x_Q - \delta \sqrt{-D}} + 4 \, x^2 \\
    & = \frac{4 \, x x_Q^3 + 4 \, x D x_Q}{(x_Q - \delta \sqrt{-D})^2} - 4 \, x \frac{x_Q^2 - D}{x_Q - \delta \sqrt{-D}} + 4 \, x^2 \\
    & = 4 \, x \frac{x_Q^2 + \delta \sqrt{-D} x_Q - x_Q^2 + D}{x_Q - \sqrt{-D}} + 4 \, x^2 \\
    & = 4 \, x^2 - \delta \sqrt{-D} x
    = \left( \frac{y}{\sqrt{\gamma}} \right)^2.
  \end{array}
  \end{equation*}
  So for the right choice of~${ \sqrt{\gamma} }$, we
  have~${ \beta + x_Q = \omega }$.
  
\end{proof}

Now we return to our premise. Suppose we have two
points~${ P, Q \in E_D(\Q) }$ which have the same class
modulo~${ [2] E_D(\Q) }$, but~${ a_P \ne a_Q }$. Note that (in $\Q^* / (\Q^*)^2$)
\begin{equation*}
  [a_P] = [A_P] = [x_P] = [y_p^2x_P] = [x_P^2 + D]
  = \left\{
    \begin{array}{ll}
      [\alpha_{+1}(P)] [\alpha_{-1}(P)] & \text{ if } -D \text{ is a square} \\
      {[\operatorname{Norm}(\alpha_{+1}(P))]} & \text{ otherwise},
    \end{array}
  \right.
\end{equation*}
and similarly for~${ [a_Q] }$.
Therefore, by Proposition~\ref{thm:descent}, we get ${ [a_P] = [a_Q] \in \Q^* / (\Q^*)^2 }$.
Without loss of generality, we may
then assume that there is a prime number~${ p }$ such that
\begin{equation*}
  0 \le \ord_p A_P < \min\{\ord_p A_Q, \ord_p D\}.
\end{equation*}
The fact that~${ P }$ and~${ Q }$ have the same class
modulo~${ [2] E_D(\Q) }$ implies
that~${ \alpha_{+1}( P ) = \alpha_{+1}( Q ) }$.
Note that~${ \alpha_{+1}( P ) }$ is the class of~${ x_P + \sqrt{-D} }$,
which is the same class as~${ A_P + B_P^2 \sqrt{-D} }$, and
similarly~${ \alpha_{+1}(Q) }$ is the class
of~${ A_Q + B_Q^2 \sqrt{-D} }$. The classes being the same implies
that
\begin{equation*}
  \left( A_P + B_P^2 \sqrt{-D} \right) \left( A_Q + B_Q^2 \sqrt{-D} \right)
  = \left( b + c \sqrt{-D} \right)^2
\end{equation*}
for some~${ b, c \in \Q }$. Since the left hand side is actually in
the ring of integers of~${ \Q(\sqrt{-D}) }$, we can in fact say
that~${ b, c \in \frac{1}{2} \Z }$
with~${ b \not\in \Z \iff c \not\in\Z }$. Writing out the equation in
the case that~${ -D }$ is not a square in~${ \Z }$, gives us
\begin{equation*}
  \begin{array}{rl}
    A_P B_Q^2 + A_Q B_P^2 & = 2 \, b c \text{ and} \\
    A_P A_Q - D B_P^2 B_Q^2 & = b^2 - D c^2.
  \end{array}
\end{equation*}
Since the left hand side of the first equation is an integer, we find
that~${ b, c \in \Z }$.
Furthermore, as~${ \ord_p A_Q > \ord_p A_P \ge 0 }$, we
have~${ \ord_p B_Q = 0 }$ and the first equation tells us that
\begin{equation*}
  \ord_p A_P = \ord_p 2 + \ord_p b + \ord_p c.
\end{equation*}
This implies that~${ \ord_p b \ge \frac{1}{2} \ord_p D > 0 }$, as
otherwise~${ b^2 }$ would be the only term in the second equation with
the smallest order of~${ p }$. Note
that~${ [ a_P ] = [ a_Q ] \in \Q^* / (\Q^*)^2 }$ implies
that~${ \min\{\ord_p A_Q, \ord_p D\} - \ord_p A_P }$ is a multiple
of~${ 2 }$, so
\begin{equation*}
  \frac{1}{2} \ord_p D \le \ord_p 2 + \ord_p b + \ord_p c = \ord_p A_P \le \ord_p D - 2,
\end{equation*}
which would imply that~${ \ord_p D \ge 4 }$, contradicting one of our
initial assumptions.
Therefore, indeed the value of~${ a_P }$ only
depends on the class of~${ P }$ modulo~${ [2] E_D(\Q) }$ in
case~${ -D }$ is not a square.

In the case that~${ -D }$ is a square, the part
about~${ \min\{\ord_p A_Q, \ord_p D\} - \ord_p A_P }$ being a multiple
of~${ 2 }$ remains true.
Since~${ -D }$ is a square, this immediately
shows that~${ \ord_p A_Q \ge \ord_p D = 2 }$ and~${ \ord_p A_P = 0
}$. This implies that~${ \ord_p(A_P + B_P^2 \sqrt{-D}) = 0 }$
and~${ \ord_p(A_Q + B_Q^2 \sqrt{-D}) = 1 }$, but this contradicts the
fact that their product should be a square. Therefore the result also
holds when~${ -D }$ is a square.

\section{Modularity results}
\label{sec:modularity}

If $a$ is a square, then~\eqref{eqn E gamma} is defined over $\Q$ for any twisting parameter $\gamma \in \Q^*$,
hence well-known to be modular~\cite{BreuilEtAl01}.
Otherwise we need to use the more general modularity
of~${ \Q }$-curves for~\eqref{eqn E gamma}. Ribet~\cite{Ribet04} proved that
every~${ \Q }$-curve~${ E }$ without complex multiplication is a
quotient over~${ \overline{\Q} }$ of a~${ \Q }$-simple abelian
variety~${ A }$ of~${ \GL_2 }$-type, and that each such
variety~${ A }$ is isogenous to an abelian variety associated with a
newform by the Shimura construction. The latter result was proven
based on the Serre conjectures, now proven by Khare and
Wintenberger~\cite{KhareWintenberger09}.

If we denote by~${ \phi_\sigma : \prescript{\sigma}{}{E} \to E }$
for~${ \sigma \in G_\Q }$ a (locally constant) choice of isogenies for
the~${ \Q }$-curve~${ E }$ without CM, then an abelian variety~${ A }$
can be constructed from a map~${ \beta : G_\Q \to \overline{\Q}^* }$
satisfying
\begin{equation}
  \label{eqn beta c}
  \beta(\sigma) \beta(\tau) \beta(\sigma \tau)^{-1} =: c_\beta(\sigma, \tau) = c_E(\sigma, \tau) := \phi_\sigma \prescript{\sigma}{}{\phi_\tau} \phi_{\sigma \tau}^{-1} \quad \text{for all } \sigma, \tau \in G_\Q,
\end{equation}
which is called a splitting map. Quer~\cite{Quer00} shows how to use the
degree map~${ d : \sigma \mapsto \degree \phi_\sigma }$ to compute
candidates for splitting
characters~${ \varepsilon : \sigma \to \beta(\sigma)^2 d(\sigma)^{-1}
}$ associated with a splitting map~${ \beta }$.

Suppose~${ E }$ is completely defined over a Galois number
field~${ K }$, i.e.~${ E }$ and the isogenies~${ \phi_\sigma }$ are
defined over ~${ K }$. If we find a candidate~${ \beta }$ which
factors over~${ G_\Q^K }$ from a splitting character candidate as
above, it might not be true that~${ c_\beta = c_E }$. However, if we
can find an~${ \alpha : \G_\Q^K \to K^* }$ such that the coboundary
of~${ \alpha }$ is~${ c_\beta c_E^{-1} }$, then the theory
in~\cite{Quer00} can be used to find a twist~${ E^{\gamma} }$ of~${ E }$
such that~${ c_\beta = c_{E^\gamma} }$. Here~${ \gamma }$
satisfies~${ \prescript{\sigma}{}{\gamma} = \alpha(\sigma)^2 \gamma }$
for all~${ \sigma \in G_\Q^K }$. The twist~${ \gamma }$ is unique up
to multiplication by elements from~${ \Q^* (K^*)^2 }$.

Proposition~5.1 and Proposition~5.2 in~\cite{Quer00} show that
if~${ E }$ is completely defined over an abelian Galois number field~${ K }$ and a splitting
map~${ \beta }$ factors over~${ G_\Q^K }$, then the restriction of
scalars~${ \residue_\Q^K E }$ is isogenous to a product
of~${ \Q }$-simple abelian varieties of~${ \GL_2 }$-type that are
pairwise non~${ \Q }$-isogenous. Furthermore, all these abelian
varieties arise from splitting maps of the form~${ \chi \beta }$
with~${ \chi : G_\Q^K \to \overline{\Q}^* }$ a character. This makes
all the corresponding newforms twists of one another, allowing one to
compute their levels using formulas that relate the conductor
of~${ E }$, the conductor of~${ \residue_\Q^K E }$, and the levels of
twists of newforms. The character of the newforms can be obtained from
Lemma~4.3 in~\cite{Ribet04}, which shows the character is the inverse of
the associated splitting character~${ \varepsilon }$.

A newform~${ f }$ associated to the abelian variety arising from a
splitting map~${ \beta }$ has~${ \lambda }$-adic and mod~${ \lambda }$
representations~${ \rho_{\lambda}^f : G_\Q \to \GL_2(L_{\lambda}) }$
and~${ \overline{\rho}_{\lambda}^f : G_\Q \to \GL_2(\F_\lambda) }$,
where~${ \lambda }$ is a finite prime of the coefficient field~${ L }$
of~${ f }$. The Galois representations~${ \rho_{\lambda}^f }$
and~${ \overline{\rho}_{\lambda}^f }$ can be related back to the
splitting map~${ \beta }$. In particular, the field~${ L }$ is also the
number field generated by the values of~${ \beta }$, and its degree is
equal to the dimension of the associated abelian variety~${ A }$.
In~\cite[Section 2.10]{LangenThesis}, one can find details about
how the data above can be used to compute the trace of Frobenius
elements for these representations.

To apply level lowering results, we will need that the Galois
representation~${ \overline{\rho}_{\lambda}^f }$ is absolutely
irreducible.
Note that irreducible and absolutely irreducible are equivalent,
as the representation~${ \overline{\rho}_\lambda^f }$ is odd
(e.g. by Theorem~3.2 in~\cite{Ribet04}).
We can derive irreducibility from the relation with the~${ \Q }$-curve.

\begin{theorem}
  \label{thm:QcurveIrreducibility}
  Let~${ E }$ be a~${ \Q }$-curve without~${ CM }$ defined
  over~${ \Q(\sqrt{a}) }$ with~${ a \in \Z \setminus \{1\} }$
  squarefree. Suppose that the corresponding
  isogeny~${ \prescript{\sigma}{}{E} \to E }$ has degree~${ 2 }$
  (it does not have to be defined over~${ \Q(\sqrt{a}) }$).
  Let~${ f }$ be
  a newform corresponding to an abelian variety of~${ \GL_2 }$-type
  that arises from a splitting map~${ \beta }$ associated
  with~${ E }$, and let~${ \lambda }$ be a finite prime of the
  coefficient field~${ L }$ of~${ f }$ with the characteristic
  of~${ \F_\lambda }$ equal to~${ l > 3 }$.
  If~${ \overline{\rho}_\lambda^f : G_\Q \to \GL_2(\F_\lambda) }$ is
  reducible, then we have that either
  \begin{itemize}
  \item ${ l = 3 }$ and
    \begin{align*}
      j(E) = \, & 2^6 y^{-2} \left( 4\,x - 7\,y \right)^{-6} \\
      & \cdot \Big(\big( 512\,x^{8} - 6\,016\,x^{7} y + 78\,176\,x^{6} y^{2} + 987\,032\,x^{5} y^{3} \\
      & \phantom{\cdot\Big(\big(} + 30\,371\,282\,x^{4} y^{4} + 97\,063\,160\,x^{3} y^{5} \\
      & \phantom{\cdot\Big(\big(} + 226\,082\,780\,x^{2} y^{6} + 227\,965\,064\,x y^{7} \\
      & \phantom{\cdot\Big(\big(} + 291\,927\,773\,y^{8} \big) \\
      & \phantom{\cdot\Big(} + 2 \sqrt{a} \left(x - 22\,y\right) \left(x + 5\,y\right)^2 \\
      & \phantom{\cdot\Big(+}\cdot\big( 256\,x^{4} - 64\,x^{3} y +
      65\,616\,x^{2} y^{2} \\
      & \phantom{\cdot\Big(+\cdot\big(}+ 80\,372\,x y^{3} + 187\,783\,y^{4} \big)
      \Big),
    \end{align*}
    for some~${ x, y \in \Q }$ with~${ x^2 + 2\,y^2 = a }$;
  \item ${ l = 5 }$ and
    \begin{align*}
      j(E) = \, & 2^6 y^{-2} \left(4\,x - 3\,y\right)^{-10} \\
      & \cdot\Big(\big( 131\,072\,x^{12} - 1\,015\,808\,x^{11} y + 15\,802\,368\,x^{10} y^{2} \\
      & \phantom{\cdot\Big(\big(} + 303\,943\,680\,x^{9} y^{3} + 8\,502\,563\,840\,x^{8} y^{4} \\
      & \phantom{\cdot\Big(\big(} + 41\,661\,192\,832\,x^{7} y^{5} + 122\,507\,172\,512\,x^{6} y^{6} \\
      & \phantom{\cdot\Big(\big(} + 219\,682\,233\,088\,x^{5} y^{7} + 344\,561\,617\,040\,x^{4} y^{8} \\
      & \phantom{\cdot\Big(\big(} + 329\,235\,309\,720\,x^{3} y^{9} + 342\,028\,231\,098\,x^{2} y^{10} \\
      & \phantom{\cdot\Big(\big(} + 150\,869\,431\,408\,x y^{11} + 111\,226\,255\,277\,y^{12}\big) \\
      & \phantom{\cdot\Big(} + 2 \sqrt{a} \left(2\,x + 11\,y\right)^2 \left(4\,x^{3} - 84\,x^{2} y - 37\,x y^{2} - 122\,y^{3}\right) \\
      & \phantom{\cdot\Big(+} \cdot \big(4\,096\,x^{6} + 7\,168\,x^{5} y +
      1\,058\,560\,x^{4} y^{2} + 2\,349\,440\,x^{3} y^{3} \\
      & \phantom{\cdot\Big(+\big(} + 4\,841\,440\,x^{2} y^{4} + 2\,594\,668\,x
      y^{5} + 3\,767\,779\,y^{6}\big) \Big),
    \end{align*}
    for some~${ x, y \in \Q }$ with~${ x^2 + y^2 = a }$;
  \item ${ l = 7 }$
    and~${ j(E) = -3375, \frac{-10\,529 \pm 16\,471 \sqrt{-7}}{8},
      \frac{56\,437\,681 \pm 1\,875\,341 \sqrt{-7}}{32\,768} }$;
  \item ${ l = 13 }$ and~${ j(E) = 3\,448\,440\,000 \pm 956\,448\,000 \sqrt{13} }$; or
  \item ${ l = 11 }$ or~${ l > 13 }$, and~${ E }$ has potential good
    reduction at all primes of characteristic~${ > 3 }$.
  \end{itemize}
\end{theorem}
\begin{proof}
  The last case, when~${ l = 11 }$ or~${ l > 13 }$, follows directly from
  Proposition~3.2 in~\cite{Ellenberg04}.
  We sketch the proof for the remaining cases; 
  for more details we refer to~\cite[Section 2.11]{LangenThesis}.
   
  First of all we can show that~${ \overline{\rho}_\lambda^f }$ is
  isomorphic
  to~${ \overline{\rho}_l^E : G_{\Q(\sqrt{a})} \to \Aut(E[l]) \cong
    \GL_2(\F_l) }$ when restricted to~${ G_{\Q(\sqrt{a})}
  }$. Therefore~${ \overline{\rho_{\lambda}^f} }$ being reducible
  implies~${ E }$ has an~${ l }$-isogeny. Since~${ E }$ also has
  a~${ 2 }$-isogeny it must correspond to a~${ \Q(\sqrt{a}) }$-point
  on~${ X_0(2 l) }$. As the~${ 2 }$-isogeny maps~${ E }$ to its Galois
  conjugate, the Fricke involution~${ w_2 }$ on~${ X_0(2) }$ maps this
  point on~${ X_0(2 l) }$ to its Galois conjugate. Therefore~${ E }$
  corresponds to a~${ \Q }$-point on the
  curve~${ C_l := X_0(2 l) / w_2 }$.

  Using Magma~\cite{Magma} we can find defining equations
  for~${ X_0(2 l) }$ and thus for~${ C_l }$. For~${ l = 7 }$
  and~${ l = 13 }$ the curve~${ C_l }$ has genus~${ 1 }$ and only
  finitely many~${ \Q }$-points. The~${ j }$-invariants of the
  corresponding points on~${ X_0(2 l) }$ are the ones listed in the
  theorem. For~${ l = 3 }$ and~${ l = 5 }$ both~${ X_0(2 l) }$
  and~${ C_l }$ are rational curves. By parameterizing
  the~${ \Q(\sqrt{a}) }$-points on~${ X_0(2 l) }$, we get
  the~${ j }$-invariants for these cases. The additional condition
  ensures they correspond to~${ \Q }$-points on~${ C_l }$.
\end{proof}
\begin{corollary}
  \label{thm:freycurveirreducible}
  Let~${ P \in E_D(\Q) \setminus \{ \mathcal{O}, T \} }$ be a point
  such that~${ B_P \ne \pm 1 }$ is an~${ l }$-th power with~${ l }$ an
  odd prime number and~${ a_P }$ not a square. Let~${ f }$ be a
  newform corresponding to an abelian variety of~${ \GL_2 }$-type that
  arises from a splitting map for~${ E_{a_P, z_P, w_P}^\gamma }$. For any
  prime~${ \lambda \mid l }$ in the coefficient field of~${ f }$, the
  Galois
  representation~${ \overline{\rho}_\lambda^f : G_\Q \to \GL_2(\F_l)
  }$ is irreducible if either
  \begin{itemize}
  \item ${ l = 3 }$ and~${ a }$ is not the norm of an element
    in~${ \Q(\sqrt{-2}) }$;
  \item ${ l = 5 }$ and~${ a }$ is not the norm of an element
    in~${ \Q(\sqrt{-1}) }$;
  \item ${ l = 7 }$ or~${ l = 13 }$; or
  \item ${ l = 11 }$ or~${ l > 13 }$, and~${ B_P }$ is divisible by a
    prime number~${ p > 3 }$.
  \end{itemize}
\end{corollary}
\begin{proof}
  Since~${ l \ge 3 }$ and~${ B_P \ne \pm 1 }$, there must be a
  prime~${ p \mid B }$ satisfying the condition of
  Proposition~\ref{thm:potentialmultiplicative}. By
  Corollary~\ref{thm:nonCM} the curve~${ E_{a_P, z_P, w_P}^\gamma }$
  therefore does not have complex multiplication allowing the
  application of Theorem~\ref{thm:QcurveIrreducibility}. The cases
  for~${ l = 3 }$ and~${ l = 5 }$ immediately follow.
  Since~${ j_{a_P, z_P, w_P} }$ is not integral, the case~${ l = 13 }$
  is also clear. For the case~${ l = 7 }$ we still need to check
  two~${ j }$-invariants, but by solving for which points~${ P }$ we
  have that~${ j_{a_P, z_P, w_P} }$ is one of the non-integral
  $ j $-invariants listed, we find that~${ B_P }$ is not a~$ 7 $-th
  power. For the cases~${ l = 11 }$ and~${ l > 13 }$ the fact that a
  prime~${ p \mid B_P }$ exists with~${ p > 3 }$, tells us we have
  potential multiplicative reduction at that prime.
\end{proof}

Now we show that, given irreducibility, we can apply level lowering
results to the Galois representation~${ \overline{\rho}_\lambda^f }$.
\begin{proposition}
  \label{thm:levelLowering}
  Let~${ P \in E_D(\Q) \setminus \{ \mathcal{O}, T \} }$ be a point
  such that~${ B_P \ne \pm 1 }$ is an~${ l }$-th power with~${ l }$ an
  odd prime number. Suppose that~${ E = E_{a_P, z_P, w_P}^\gamma }$ is defined
  over a number field~${ K }$ and
  that~${ f \in \mathcal{S}_2(N, \chi) }$ is a newform such that
  \begin{equation*}
    \overline{\rho}_{\lambda}^f |_{G_K} \cong \overline{\rho}_l^E : G_K \to \Aut(E[l]) \cong \GL_2(\F_l)
  \end{equation*}
  for some prime~${ \lambda \mid l }$ in the coefficient field
  of~${ f }$. Let~${ S }$ be the set of prime numbers that ramify
  in~${ K }$, divide~${ 2 D }$, or divide the norm of~${ \gamma
  }$.
  Denote~${ \tilde{N} := \prod_{p \in S} p^{\ord_p N} }$.
  If~${ \rho_{\lambda}^f }$ is irreducible, then there exists a
  newform~${ g \in \mathcal{S}_2(\tilde{N}, \chi) }$ and a
  prime~${ \lambda' \mid l }$ in the coefficient field of~${ g }$ such
  that
  \begin{equation*}
    \overline{\rho}_{\lambda'}^g \cong \overline{\rho}_{\lambda}^f : G_\Q \to \GL_2(\overline{\F_l}).
  \end{equation*}
\end{proposition}
\begin{proof}
  For any prime number~${ p }$ not in~${ S }$, we know that~${ p }$
  does not ramify in~${ K }$, hence the ramification
  subgroup~${ I_p \subseteq G_\Q }$ is a subgroup of~${ G_K }$. This
  implies
  that~${ \overline{\rho}_{\lambda}^f | I_p \cong \rho_l^E | I_p
  }$. Note that~${ E }$ has good or multiplicative reduction at each
  prime~${ \mathfrak{p} }$ above~${ p }$ by
  Corollary~\ref{thm:nonadditive} and is also a minimal model at those
  primes.
  Furthermore, the assumption that~${ B_P }$ is an~${ l }$-th
  power shows
  that~${ l \mid \ord_{\mathfrak{p}} \Delta_{a_P, z_P, w_P}^\gamma }$
  as~${ p \nmid 2 D }$ and~${ p }$ also does not divide the norm
  of~${ \gamma }$. It is well known that this implies
  that~${ \rho_l^E | I_p }$ is trivial when~${ p \ne l }$ and finite
  flat if~${ p = l }$.
  Standard level lowering results now apply to obtain the newform $g$.
\end{proof}

\section{The case ${ a = 1 }$}
\label{sec:aeq1}

The case~${ a = 1 }$ is special as shown by the following result.
\begin{proposition}
  \label{thm:a=1foreven}
  For any~${ P \in [2]E_D(\Q) \setminus \{ \mathcal{O} \} }$, we
  have~${ a_P = 1 }$.
\end{proposition}
\begin{proof}
  We use the
  homomorphisms~${ \alpha_\delta : E_D(\Q) \to \Q(\sqrt{-D})^* /
    (\Q(\sqrt{-D})^*)^2 }$ from Proposition~\ref{thm:descent}. For
  any~${ P \in [2]E_D(\Q) }$ and~${ \delta \in \{ \pm 1 \} }$, we know
  that~${ \alpha_\delta(P) = [1] }$.
  Furthermore, we have seen that
  \begin{equation*}
    \Q^* / (\Q^*)^2 \ni [a_P] = \left\{
      \begin{array}{cl}
        [\alpha_{+1}(P)] [\alpha_{-1}(P)] & \text{if } -D \text{ is a square} \\
        {[\operatorname{Norm}(\alpha_{+1}(P)]} & \text{otherwise,}
      \end{array}
    \right.
  \end{equation*}
  so~${ a_P }$ is a square. We know
  that~${ \alpha_{+1}(P) = [A_P + B_P^2 \sqrt{-D}] }$,
  so~${ A_P + B_P^2 \sqrt{-D} }$ is a square
  in~${ \mathcal{O}_{\Q(\sqrt{-D})} }$. We can
  write~${ A_P + B_P^2 \sqrt{-D} = (b + c \sqrt{-D})^2 }$
  with~${ b, c \in \frac{1}{2} \Z }$ if~${ -D }$ is not a square to
  find that
  \begin{equation*}
    \left\{
      \begin{array}{l}
        A_P = b^2 - D c^2 \\
        B_P^2 = 2 b c.
      \end{array}
    \right.
  \end{equation*}
  Since~${ b + c \sqrt{-D} \in \mathcal{O}_{\Q(\sqrt{-D})} }$, we may
  assume from the second equation that~${ b, c \in \Z }$. Now any
  prime dividing~${ a_P }$ divides both~${ A_P }$ and~${ D }$, hence
  divides also~${ b }$ and thus~${ B_P }$.
  Since~${ \gcd(A_P, B_P) = 1 }$, this shows us that~${ a_P }$ = 1 in
  this case.

  If~${ -D }$ is a square, we know
  that~${ A_P + B_P^2 \sqrt{-D} = b^2 }$ for some~${ b \in \Z }$. Note
  that for any prime~${ p \mid a_P }$, we have~${ p \nmid B_P }$ and
  therefore
  \begin{equation*}
    \ord_p A_P \ge \ord_p a_P \ge 2 > \frac{1}{2} \ord_p D = \ord_p (B_P^2 \sqrt{-D})
  \end{equation*}
  by our assumptions on~${ D }$. This implies that
  \begin{equation*}
    2 \ord_p b = \ord_p (A_P + B_P^2 \sqrt{-D}) = \frac{1}{2} \ord_p D > 0,
  \end{equation*}
  which contradicts the previous equation. Therefore~${ a_P = 1 }$ also
  in this case.
\end{proof}
This result shows that for any elliptic divisibility sequence
generated by~${ P }$, the curve~${ E_{1, z, w}^\gamma }$ is the Frey
curve corresponding to the points~${ P_{2 m}, m \in \Z }$. Note that
the curve~${ E_{1, z, w}^\gamma }$ is an elliptic curve defined
over~${ \Q }$ when~${ \gamma \in \Q^* }$ and does not depend
on~${ D }$. We shall prove some general results about this curve here.

First we will try to find the twist~${ \gamma \in \Q^* }$ such that
the conductor of~${ E_{1, z, w}^\gamma }$ is as small as possible.
For this, we may assume that~${ \gamma }$ is a squarefree integer. Note
that if an odd prime~${ p }$ divides~${ \gamma }$,
then~${ c_{4, 1, z, w}^\gamma }$ and~${ \Delta_{1, z, w}^\gamma }$ are
both divisible by~${ p }$, but the model is still minimal as
either~${ p^4 \nmid c_{4, 1, z, w}^\gamma }$
or~${ p^6 \nmid \Delta_{1, z, w}^\gamma }$. So for an odd
prime~${ p }$, the conductor exponent at~${ p }$ would be at
least~${ 2 }$ if~${ p \mid \gamma }$, whereas it would be~${ 0 }$
or~${ 1 }$ otherwise by Corollary~\ref{thm:nonadditive}. We therefore
want~${ \gamma \in \{ 1, -1, 2, -2 \} }$.

Table~\ref{tab:conductorat2} gives the possible conductor exponents
at~${ 2 }$ for the curves~${ E_{1, z, w}^\gamma }$
with~${ \gamma }$ in~${\{ 1, -1, 2, -2 \} }$. To compute these conductor
exponents, we used that the curve~${ E_{1, z, w}^\gamma }$
is~${ 2 }$-isogenous over~${ \Q }$ to the curve
\begin{equation*}
  \tilde{E}_{1, z, w}^{\gamma} : Y^2 = X^3 - 8 \, z \gamma X^2+ 8\left(z^2 - w \right) \gamma^2 X.
\end{equation*}
Since isogenous curves have isomorphic~${ l }$-adic Galois
representations, their conductors are by definition the same. In the
framework~\cite{framework} we can use this to perform the conductor
computation on both curves simultaneously, and use the result from the
one that finishes first for each pair~${ z, w }$.

\begin{table}[h]
  \makebox[\textwidth][c]{
  \begin{tabular}{c|c|c|c|l}
    ${ \gamma = 1 }$ & ${ \gamma = -1 }$ & ${ \gamma = 2 }$ & ${ \gamma = -2 }$ & \\
    \hline
    8 & 8 & 8 & 8 & if~${ \ord_2 (w^2 - z^4) = 0 }$ \\
    7 & 7 & 7 & 7 & if~${ \ord_2 (w^2 - z^4) = 3 }$ \\
    4 & 3 & 6 & 6 & if~${ w + z^2 \equiv 8 }$ (mod $32$) and $z \equiv 1$ (mod $4$) \\
    3 & 4 & 6 & 6 & if~${ w + z^2 \equiv 8 }$ (mod $32$) and $z \equiv 3$ (mod $4$) \\
    2 & 4 & 6 & 6 & if~${ w + z^2 \equiv 24 }$ (mod $32$) and $z \equiv 1$ (mod $4$) \\
    4 & 2 & 6 & 6 & if~${ w + z^2 \equiv 24 }$ (mod $32$) and $z \equiv 3$ (mod $4$) \\
    6 & 6 & 4 & 2 & if~${ w - z^2 \equiv 8 }$ (mod $32$) and $z \equiv 1$ (mod $4$) \\
    6 & 6 & 2 & 4 & if~${ w - z^2 \equiv 8 }$ (mod $32$) and $z \equiv 3$ (mod $4$) \\
    6 & 6 & 3 & 4 & if~${ w - z^2 \equiv 24 }$ (mod $32$) and $z \equiv 1$ (mod $4$) \\
    6 & 6 & 4 & 3 & if~${ w - z^2 \equiv 24 }$ (mod $32$) and $z \equiv 3$ (mod $4$) \\
    5 & 5 & 6 & 6 & if~${ \ord_2 (w + z^2) = 4 }$ \\
    6 & 6 & 5 & 5 & if~${ \ord_2 (w - z^2) = 4 }$ \\
    3 & 4 & 6 & 6 & if~${ \ord_2 (w + z^2) = 5, 6 }$ and $z \equiv 1$ (mod $4$) \\
    4 & 3 & 6 & 6 & if~${ \ord_2 (w + z^2) = 5, 6 }$ and $z \equiv 3$ (mod $4$) \\
    6 & 6 & 4 & 3 & if~${ \ord_2 (w - z^2) = 5, 6 }$ and $z \equiv 1$ (mod $4$) \\
    6 & 6 & 3 & 4 & if~${ \ord_2 (w - z^2) = 5, 6 }$ and $z \equiv 3$ (mod $4$) \\
    0 & 4 & 6 & 6 & if~${ \ord_2 (w + z^2) = 7 }$ and $z \equiv 1$ (mod $4$) \\
    4 & 0 & 6 & 6 & if~${ \ord_2 (w + z^2) = 7 }$ and $z \equiv 3$ (mod $4$) \\
    6 & 6 & 4 & 0 & if~${ \ord_2 (w - z^2) = 7 }$ and $z \equiv 1$ (mod $4$) \\
    6 & 6 & 0 & 4 & if~${ \ord_2 (w - z^2) = 7 }$ and $z \equiv 3$ (mod $4$) \\
    1 & 4 & 6 & 6 & if~${ \ord_2 (w + z^2) \ge 8 }$ and $z \equiv 1$ (mod $4$) \\
    4 & 1 & 6 & 6 & if~${ \ord_2 (w + z^2) \ge 8 }$ and $z \equiv 3$ (mod $4$) \\
    6 & 6 & 4 & 1 & if~${ \ord_2 (w - z^2) \ge 8 }$ and $z \equiv 1$ (mod $4$) \\
    6 & 6 & 1 & 4 & if~${ \ord_2 (w - z^2) \ge 8 }$ and $z \equiv 3$ (mod $4$)
  \end{tabular}}
  \caption{The conductor exponent of~${ E_{1, z, w}^\gamma }$ at~$2$ for various values of~${ \gamma }$.}
  \label{tab:conductorat2}
\end{table}

\begin{remark}
  \label{rem:unneeded exps}
  Note that the only assumption used to generate
  Table~\ref{tab:conductorat2} is that~${ 2 }$ does not divide
  both~${ z }$ and~${ w }$. When~${ z }$ and~${ w }$ correspond to a
  point on~${ E_D }$, we have that
  \begin{equation*}
    (w + z^2) (w - z^2) = w^2 - z^4 = D B^4,
  \end{equation*}
  so we obtain further limitation based on the order of~${ 2 }$
  in~${ D }$ and~${ B }$. In particular, if we assume~${ B }$ is
  an~${ l }$-th power with~${ l > 1 }$, the cases with conductor
  exponents~${ 2 }$,~${ 3 }$, and~${ 5 }$ do not occur.
\end{remark}
We will assume that whenever we use~${ E_{1, z, w}^\gamma }$ from now
on,~${ \gamma }$ is chosen such that the conductor exponent at~${ 2 }$
is as small as possible. Note that to do this, ${ \gamma }$ might
depend on the value of~${ z }$ and~${ w }$ modulo~${ 2^8 }$.

We will also show here that the Galois representations associated to
the case~${ a = 1 }$ are irreducible.
\begin{theorem}
  \label{thm:am=1Irreducible}
  For all prime numbers~${ l > 5 }$ and~${ \gamma \in \Q^* }$, the
  mod~${ l }$ Galois
  representation~${ \rho_l^{E_{1, z, w}^\gamma} : G_\Q \to \Aut(E_{1, z, w}^\gamma [ l ]) \cong \GL_2(\F_l) }$
  is irreducible.
\end{theorem}

\begin{proof}
  Note that~${ E_{1, z, w} }$ has a rational~$ 2 $-torsion point.
  If~${ \rho_l^{E_{1, z, w}^\gamma} }$ would be reducible, then $E$
  would correspond to a non-cuspidal~${ \Q }$-rational point
  on~${ X_0(2 l) }$. The well known deep fact that $X_0(2l)$ does not
  have such points for~${ l > 7 }$ (see e.g. \cite[Theorem
  22-(ii)]{Dahmen08} for precise references) immediately proves the
  result for~${ l > 7 }$.
  
  For~${ l = 7 }$, we will check if
  the~${ j }$-invariants of points on~${ X_0(2 l) }$ can
  match~${ j_{1, z, w} }$. Note that we have
  \begin{equation*}
    j_{1, z, w} = 2^6 \frac{(3 \, t - 5)^3}{(t - 1) (t + 1)^2},
  \end{equation*}
  after rewriting with~${ t = \frac{w}{z^2} }$.
  The previous reference basically also notes that the image to the $j$-line
  of the~${ \Q }$-rational points on~${ X_0(14) }$
  equals~${ \{\infty, -3375, 16\,581\,375\} }$.
  Equating these $j$-values to~${ j_{1,z,w} }$, we
  find~${ t = \pm 1, \pm \frac{65}{63} }$. Since the denominator
  of~${ t }$ is a square, the only option is~${ t = \pm 1 }$, which
  would give~${ z = \pm 1 }$ and~${ w = \pm 1 }$, meaning~${ B = 0
  }$. This does not correspond to a
  point~${ P \in E_D(\Q) \setminus \{ \mathcal{O}, T \} }$.
\end{proof}

We can extend Theorem~\ref{thm:am=1Irreducible} to~${ l = 3, 5 }$ for
specific~${ D }$ as follows. First assume that the mod~${ l }$ Galois
representation~${ \rho_l^{E_{1, z, w}^\gamma} : G_\Q \to \Aut(E_{1, z,
    w}^\gamma [ l ] })$ is reducible, then~${ E_{1, z, w}^\gamma }$
corresponds to a~${ \Q }$-point on~${ X_0(l) }$. Note
that~${ X_0(l) }$ is a rational curve for~${ l = 3, 5 }$, so we can
parameterize it with a single parameter~${ s }$. We
let~${ j_l : X_0(l) \to X(1) }$ be the~${ j }$-invariant, in which case
we get that the curve~${ E_{1, z, w}^\gamma }$ should correspond to a
point on the curve
\begin{equation*}
  C_l : j_{1, z, w}(t) = j_l(s).
\end{equation*}
Computing this curve explicitly in Magma~\cite{Magma}, we see that it
also is a rational curve, hence it has a
parameterization. Parameterizing with respect to
points~${ [ x : y ] }$ on~${ \Proj^1 }$ we get that
\begin{equation*}
  \frac{w}{z^2} = t = \frac{w_l(x, y)}{z_l(x, y)}.
\end{equation*}
Since any point of~${ \Proj^1 }$ can be written with coprime integer
coordinates, we thus know there are coprime~${ a, b \in \Z }$ and
some~${ c \in \Q^* }$ such that
\begin{equation*}
  \left\{
    \begin{array}{ll}
      c w & = w_l(a, b), \\
      c z^2 & = z_l(a, b).
    \end{array}
  \right.
\end{equation*}
By rescaling we may assume~${ w_l }$ and~${ z_l }$ have integer
coefficients, in which case the denominator of~${ c }$ must be one
as~${ w }$ and~${ z^2 }$ are coprime. Therefore, we may also
assume~${ c \in \Z \setminus \{ 0 \} }$.

From these equations, it follows that
\begin{equation*}
  c^2 D B^4 = (c w)^2 - (c z^2)^2 = w_l(a, b)^2 - z_l(a, b)^2.
\end{equation*}
Computing the right hand side in Magma, we actually see that
\begin{equation}
  \label{eq:BprodOfLinearFactors}
  c^2 D B^4 = 2^8 c_1(a, b)^l c_2(a, b)^l c_3(a, b) c_4(a, b)
\end{equation}
with~${ c_1, c_2, c_3, c_4 }$ all linear factors. By choosing a
different parameterization of~${ \Proj^1 }$, we may assume that
\begin{equation*}
  \begin{array}{ll}
    c_1(a, b) & = a, \\
    c_2(a, b) & = b, \\
    c_3(a, b) & = a - b.
  \end{array}
\end{equation*}
Making this parameterization explicit in Magma, we find that
\begin{equation*}
  c_4(a, b) = \left\{
    \begin{array}{ll}
      a + 8 \, b & \text{ if } l = 3 \\
      a + 4 \, b & \text{ if } l = 5.
    \end{array}
  \right.
\end{equation*}
Since~${ a }$ and~${ b }$ are coprime, we can easily see that
the~${ c_i }$ are pairwise coprime outside primes
dividing~${ 2 \, l }$. Therefore,
equation~\eqref{eq:BprodOfLinearFactors} implies that all
the~${ c_i(a, b) }$ must be fourth powers up to factors consisting
of primes dividing~${ 2 \, l c^2 D }$, i.e. we can write
\begin{equation*}
  c_i(a, b) = a_i b_i^4 \quad \text{with } a_i, b_i \in \Z.
\end{equation*}
Note that as~${ c }$ must divide the resultant of~${ z_l }$
and~${ w_l }$, there is only a finite list of
possible~${ (a_1, a_2, a_3, a_4) }$. We can compute this list
explicitly for a given~${ D }$ using Magma. For each choice of
the~${ a_i }$, the linear relations between the~${ c_i }$ give us four
generalized Fermat equations of signature~${ (4, 4, 4) }$ that should
be satisfied. Reaching a contradiction now follows by eliminating one
such equation for each choice of~${ a_i }$.

\begin{table}[b!]
  \centering
  \begin{tabular}{c|c|c}
    $ D $ & $ l = 3 $ & $ l = 5 $ \\
    \hline
    $ -2 $ & $ \mathcal{O} $ & $ \mathcal{O} $ \\
    $ 3 $ & $ \mathcal{O} $ & $ \mathcal{O} $, $ (1, \pm 2) $, $ \left(\frac{121}{9}, \pm\frac{1\,342}{27}\right) $ \\
    $ -17 $ & $ \mathcal{O} $ & $ \mathcal{O} $ \\
    $ 125 $ & $ \mathcal{O} $, $ \left(\frac{121}{4}, \frac{1\,419}{8}\right) $ & $ \mathcal{O} $
  \end{tabular}
  \caption{Points on~${ E_D(\Q) }$ such that the mod~${ l }$ Galois representation of~${ E_{1, z, w}^\gamma }$ is irreducible if it does not correspond to that point.}
  \label{tab:l35irreducible}
\end{table}

To eliminate Fermat equations, we first check if they have local
solutions over~${ \Q_2 }$, ${ \Q_3 }$, or~${ \Q_5 }$. If all four
corresponding to a choice of the~${ a_i }$ do, we look at the quotient
\begin{equation*}
  \begin{array}{rl}
    \Proj^2 \supseteq \left\{ A X^4 + B Y^4 + C Z^4 = 0 \right\}
    & \to \left\{ Y^2 Z = X^3 + \frac{B C}{A^2} X Z^2 \right\} \subseteq \Proj^2 \\
    {[ x : y : z ]} & \mapsto [ B y^2 z : B x^2 y : -A z^3 ]
  \end{array}
\end{equation*}
and the other genus~${ 1 }$ quotients obtained by
interchanging~${ X }$,~${ Y }$, and~${ Z }$. If one of these elliptic
curves has rank~${ 0 }$, we can determine the solutions of the
corresponding Fermat equation that map to its torsion points. This in
turn gives us the possible values of~${ a }$ and~${ b }$ and hence the
values of~${ z, w }$ for which the mod~${ l }$ representation may
still be reducible. If we can do this for all remaining choices
of~${ a_i }$, we thus get an explicit list of points on~${ E_D(\Q) }$
such that if~${ E_{1, z, w}^\gamma }$ does not correspond to one of
these points, then its mod~${ l }$ Galois representation is
irreducible.

We did the above computation for some~${ D }$ that will be used in the
examples below. The points for which we could not show that the
mod~${ l }$ Galois representation is irreducible in this way, are
listed in Table~\ref{tab:l35irreducible}. The code written to do so
can easily be reused to compute this for further values of~${ D }$.

\section{Explicit examples}
\label{sec:EDS examples}

In this section we will make some choices of integers~${ D }$ and
points~${ P_1 \in E_D(\Q) }$ for which we can prove the non-existence
of~${ l }$-th powers among the~${ B_m }$ with an explicit lower bound
on the prime number~${ l }$.

\subsection{Example for ${D = 125}$}
\label{subsec:example 125}
Take~${ D = 125 }$, then~${ E_D(\Q) }$ has rank~1 and~${ T = (0, 0) }$
is the only non-trivial torsion point. Using SageMath we can compute
that~${ E_D(\Q) }$ is generated
by~${ P = \left(\frac{121}{4}, \frac{1419}{8}\right) }$ and~${ T
}$. We have~${ a_P = 1 }$, so all non-zero multiples of~${ P }$
correspond to the Frey curve~${ E_{1, z, w}^\gamma }$, where we
choose~${ \gamma }$ as in the previous section.

The level of the newforms that remain after level lowering is the
product of the~${ 2 }$-part and~${ 5 }$-part of the conductor
of~${ E_{1, z, w}^\gamma }$. The~${ 2 }$-part can be read of from
Table~\ref{tab:conductorat2}, where we note that
Remark~\ref{rem:unneeded exps} applies and~${ 2 \mid B_{m P} }$ for
all~${ m \in \Z \setminus \{0\} }$ by~\eqref{eqn:edsp}. We can compute
that the conductor exponent at~${ 5 }$ is always~${ 1 }$, so the
levels of the newforms after level lowering must be
\begin{equation*}
  \left\{
    \begin{array}{cl}
      5 & \text{if } w^2 - z^4 \equiv 2^8 \text{ (mod } 2^9) \\
      10 & \text{if } w^2 - z^4 \equiv 0 \text{ (mod } 2^9).
    \end{array}
  \right.
\end{equation*}
Note that there are no newforms of level~${ 5 }$ or level~${ 10 }$,
hence none of the~${ B_{m P} }$ with~${ m \in \Z \setminus \{0\} }$
can be~${ l }$-th powers for~${ l }$ an odd prime number.

There is one other Frey curve associated with~${ D = 125 }$, as for
any~${ m \in \Z \setminus \{ 0 \} }$ we have~${ a_{m P + T} = 125
}$. We will find the correct twist of the corresponding
Frey~${ \Q }$-curve~${ E_{125, z, w}^\gamma }$ by
studying~${ E_{125, z, w} }$ first. Note that~${ E_{125, z, w} }$ is
completely defined over~${ \Q( \sqrt{2}, \sqrt{5} ) }$, but a
splitting character is given by a character of conductor~${ 20 }$ and
order~${ 4 }$. Therefore, a complete definition field over which also a
splitting map is defined
is~${ K = \Q(\zeta_{40} + \zeta_{40}^{-1}) }$, the totally real
subfield of the cyclotomic field~${ \Q(\zeta_{40}) }$. We
compute~${ c_{E_{125, z, w}} }$ and~${ c_\beta }$ on~${ G_\Q^{K} }$.
We can use the
isomorphism~${ \left( \Z / 40 \Z \right)^* \cong G_\Q^{\Q(\zeta_{40})}
}$ to identify~${ G_\Q^K }$
with~${ \left( \Z / 40 \Z \right)^* / \{ \pm 1 \} }$, which we will
use to denote elements of~${ G_\Q^K }$ by representatives
from~${ \left( \Z / 40 \Z \right)^* }$.
\begin{center}
\begin{tabular}{c|cccccccc}
  $c_{E_{125, z, w}}$ & $\pm 1$ & $\pm 3$ & $\pm 7$ & $\pm 9$ & $\pm 11$ & $\pm 13$ & $\pm 17$ & $\pm 19$ \\
  \hline
  $\pm 1$ & $1$ & $1$ & $1$ & $1$ & $1$ & $1$ & $1$ & $1$ \\
  $\pm 3$ & $1$ & $-2$ & $-2$ & $1$ & $1$ & $-2$ & $-2$ & $1$ \\
  $\pm 7$ & $1$ & $2$ & $2$ & $1$ & $1$ & $2$ & $2$ & $1$ \\
  $\pm 9$ & $1$ & $1$ & $1$ & $1$ & $1$ & $1$ & $1$ & $1$ \\
  $\pm 11$ & $1$ & $-1$ & $-1$ & $1$ & $1$ & $-1$ & $-1$ & $1$ \\
  $\pm 13$ & $1$ & $-2$ & $-2$ & $1$ & $1$ & $-2$ & $-2$ & $1$ \\
  $\pm 17$ & $1$ & $2$ & $2$ & $1$ & $1$ & $2$ & $2$ & $1$ \\
  $\pm 19$ & $1$ & $-1$ & $-1$ & $1$ & $1$ & $-1$ & $-1$ & $1$ \\
\end{tabular}
\end{center}
\begin{center}
\begin{tabular}{c|cccccccc}
  $c_\beta$ & $\pm 1$ & $\pm 3$ & $\pm 7$ & $\pm 9$ & $\pm 11$ & $\pm 13$ & $\pm 17$ & $\pm 19$ \\ \hline
  $\pm 1$ & $1$ & $1$ & $1$ & $1$ & $1$ & $1$ & $1$ & $1$ \\
  $\pm 3$ & $1$ & $-2$ & $2$ & $1$ & $1$ & $2$ & $-2$ & $1$ \\
  $\pm 7$ & $1$ & $2$ & $2$ & $-1$ & $-1$ & $2$ & $2$ & $1$ \\
  $\pm 9$ & $1$ & $1$ & $-1$ & $-1$ & $-1$ & $-1$ & $1$ & $1$ \\
  $\pm 11$ & $1$ & $1$ & $-1$ & $-1$ & $-1$ & $-1$ & $1$ & $1$ \\
  $\pm 13$ & $1$ & $2$ & $2$ & $-1$ & $-1$ & $2$ & $2$ & $1$ \\
  $\pm 17$ & $1$ & $-2$ & $2$ & $1$ & $1$ & $2$ & $-2$ & $1$ \\
  $\pm 19$ & $1$ & $1$ & $1$ & $1$ & $1$ & $1$ & $1$ & $1$ \\
\end{tabular}
\end{center}

A simple computation shows us that the
map~${ \alpha : G_\Q^K \to \mathcal{O}_K^* }$ given by the table below
has coboundary~${ c_{E_{125, z, w}} c_\beta^{-1} }$.
\begin{center}
  \renewcommand{\arraystretch}{1.1}  
  \begin{tabular}{c||c|c|c|c}
    $\sigma$ & $\pm 1$, $\pm 19$ & $\pm 3$, $\pm 17$ & $\pm 7$, $\pm 13$ & $\pm 9$, $\pm 11$ \\
    \hline
    $\alpha(\sigma)$ & $1$ & $\zeta_{40}^{17} + \zeta_{40}^{-17}$ & $(\zeta_{40} + \zeta_{40}^{-1})^{-1}$ & $(\zeta_{40}^{3} + \zeta_{40}^{-3}) (\zeta_{40}^9 + \zeta_{40}^{-9})$
  \end{tabular}
\end{center}
It is easy to verify
that~${ \gamma = (\zeta_{40}^{1} + \zeta_{40}^{-1}) (\zeta_{40}^{2} +
  \zeta_{40}^{-2}) (\zeta_{40}^{3} + \zeta_{40}^{-3}) }$
satisfies~${ \prescript{\sigma}{}{\gamma} = \alpha(\sigma)^2 \gamma }$
for all~${ \sigma \in G_\Q^K }$, hence~${ \gamma }$ is the sought
twist. Note that~${ \gamma }$ and thereby the
twist~${ E_{125, z, w}^\gamma }$ are defined over the
field~${ K_0 = \Q(\zeta_{40}^2 + \zeta_{40}^{-2}) }$, which is the
totally real subfield of~${ \Q(\zeta_{20}) }$. Some further checking
shows that~${ E_{125, z, w}^\gamma }$ is actually completely defined
over this field, and the splitting map
for~${ c_{E_{125, z, w}^\gamma} }$ factors over~${ G_\Q^{K_0} }$ as
well.

We compute that the conductor of~${ E_{125, z, w}^{\gamma} }$ is equal
to
\begin{equation*}
  (64) \Rad_{10} \Delta_{125, z, w}^\gamma.
\end{equation*}
By computing all the splitting
maps~${ \beta: G_\Q^{K_0} \to \overline{\Q}^* }$, we find that the
restriction of scalars~${ \residue_\Q^{K_0} E_{125, z, w}^\gamma }$
must be isogenous to a product of two abelian varieties
of~${ \GL_2 }$-type. The levels of the associated newforms should be
\begin{equation*}
  \{ 1280 \Rad_{10}(w^2 - 125 z^4), 6400 \Rad_{10}(w^2 - 125 z^4) \},
\end{equation*}
where we can not a priori determine which level corresponds to which
factor. After level lowering, we thus obtain newforms of
levels~${ 1280 }$ and~${ 6400 }$ which must be twists of one
another. The characters of these newforms should be a character of
conductor~${ 20 }$ and order~${ 4 }$.

The irreducibility of the mod~${ l }$ representations needed to
perform level lowering is obtained from
Corollary~\ref{thm:freycurveirreducible}. For the
cases~${ l = 3, 7, 13 }$ this is immediately clear. For the
cases~${ l = 11 }$,~${ l > 13 }$ we note that the only
points~${ Q \in E_D(\Q) \setminus \{ \mathcal{O}, T \} }$
with~${ B_Q }$ not divisible by a prime number~${ p > 3 }$
are~${ \pm P }$. For the case~${ l = 5 }$, we do not get irreducibility
as~${ 125 }$ is the norm of~${ 11 + 2 \sqrt{-1} }$.

\begin{remark}
  One could try to use the~${ j }$-invariant~${ j_5(x, y) }$ presented
  in Theorem~\ref{thm:QcurveIrreducibility} to prove irreducibility of
  the mod~${ l }$ representation when~${ l = 5 }$. For this, one would
  parameterize the~${ x, y \in \Q }$ with~${ x^2 + y^2 = 125 }$ and
  find the solutions to~${ j_{125, z, w} = j_5(x, y) }$. One will
  however find that there are infinitely many solutions to this
  equation.
\end{remark}

When we compute the newforms of level~${ 1280 }$ and their twists of
level~${ 6400 }$, we end up with~${ 144 }$ newforms. After comparing
traces of Frobenius for the primes~${ p < 50 }$ with~${ p \ne 2, 5 }$,
we find that all but~${ 24 }$ newforms can be a priori eliminated for
all primes~${ l > 17 }$. All of the remaining newforms have
coefficient field~${ \Q(\sqrt{-1}) }$ and most likely correspond to
(pseudo) solutions~${ (z, w) }$. At the very least, there will be a
newform corresponding to the curve~${ E_{125, 0, 1}^\gamma }$, which
corresponds to the point~${ T }$, so we can not eliminate all newforms
without further assumptions.

If we look at the points~${ m P + T }$ with~${ m }$ an odd integer, we
note they are all the odd multiples
of~${ P + T = \left(\frac{500}{121}, -\frac{32\,250}{1\,331}\right)
}$. By~\eqref{eqn:edsp}, we have that~${ 11 \mid B_{m P + T} }$ for all
odd~${ m }$. If we use this information when computing the traces of
Frobenius at~${ 11 }$, we can actually eliminate all newforms
for~${ l > 3 }$ and~${ l \ne 11 }$.

Combining these results we have thus proven the following result.
\begin{theorem}
  \label{thm:D125}
  The sequence~${ B_m, m \in \Z_{>0} }$ contains no~${ l }$-th powers
  with~${ l }$ a prime number when
  \begin{itemize}
  \item ${ P_1 = \left(\frac{121}{4}, \pm\frac{1419}{8}\right) \in E_{125}(\Q) }$ and~${ l > 2 }$; or
  \item
    ${ P_1 = \left(\frac{500}{121}, \pm\frac{32\,250}{1\,331}\right) \in
      E_{125}(\Q) }$ and~${ l > 5 }$ with~${ l \ne 11 }$.
  \end{itemize}
\end{theorem}
\subsection{Example for ${D = -17}$}
Take~${ D = -17 }$, then~${ E_D(\Q) }$ has rank~2
with~${ T = (0, 0) }$ as the only non-trivial torsion point. With
SageMath~\cite{sagemath}, we can compute that there is a choice of
generators
\begin{equation*}
  \left\{
    \begin{array}{ll}
      P = & (-4, 2) \\
      Q = & (-1, 4) \\
      T = & \left(0, 0\right),
    \end{array}
  \right.
\end{equation*}
for which~${ a_P = a_Q = -1 }$.

We will first study the Frey curve~${ E_{1, z, w} }$, which
corresponds to the points~${ m P + n Q }$ with~${ m, n \in \Z }$
and~${ m + n }$ even. We use the twist~${ E_{1, z, w}^\gamma }$ from
Section~\ref{sec:aeq1} to get the lowest conductor. The level of the
newforms after level lowering will consist of the~${ 2 }$-part
and~${ 17 }$-part of this conductor. The~${ 2 }$-part can be read of
from Table~\ref{tab:conductorat2}, where we apply
Remark~\ref{rem:unneeded exps}. Since the conductor exponent
at~${ 17 }$ is always~${ 1 }$, we find that the level after level
lowering is
\begin{equation*}
  \left\{
    \begin{array}{cl}
      2^8 \cdot 17 & \text{if } \ord_2(w^2 - z^4) = 0 \\
      17 & \text{if } \ord_2(w^2 - z^4) = 8 \\
      2 \cdot 17 & \text{if } \ord_2(w^2 - z^4) > 8.
    \end{array}
  \right.
\end{equation*}
We compute the newforms of these levels and compare traces of
Frobenius for all prime numbers~${ p < 50 }$ that do not
divide~${ 2 \cdot 17 }$. This eliminates 25 of the 33 newforms of
level~${ 2^8 \cdot 17 }$ for all primes~${ l > 5 }$, but none of
level~${ 17 }$ or~${ 2 \cdot 17 }$. The non-eliminated newforms are
all rational and the corresponding elliptic curves are geometrically
isomorphic to the elliptic curves~${ E_{1, z, w}^\gamma }$
corresponding to the points~${ P - Q }$ and~${ Q - P }$, and the
pseudo
solutions~${ (z, w) = (\pm 12, \pm 145), (\pm 15, \pm 353), (\pm 23,
  \pm 495) }$. Therefore, they can not be eliminated without additional
information.

Now look at the non-zero multiples~${ P_m }$ of the point
\begin{equation*}
  P_1 = 2 P + 2 Q = \left(\frac{3\,568\,321}{451\,584},
    \frac{5\,750\,178\,337}{303\,464\,448}\right).
\end{equation*}
Note that~${ B_1 }$ is divisible by the primes~${ 2 }$,~${ 3 }$,
and~${ 7 }$, so by~\eqref{eqn:edsp}, all~${ B_m }$ with~${ m \in \Z }$
non-zero are divisible by~${ 2 }$,~${ 3 }$, and~${ 7 }$.
The fact that~${ 2 \mid B_m }$ immediately rules out all
newforms of level~${ 2^8 \cdot 17 }$. Comparing traces of Frobenius
at~${ 3 }$ and~${ 7 }$, again using the restriction~${ 3, 7 \mid B_m }$,
we can now eliminates all newforms for~${ l > 3 }$.

Next, we study the Frey~${ \Q }$-curve~${ E_{-17, z,w}^{\gamma} }$ (see
equation~\eqref{eqn E gamma}) corresponding to
points~${ m P + n Q + T }$ with~${ m + n }$ even that have
an~${ l }$-th power in their denominator.

We start with~${ \gamma = 1 }$ and find the correct twist later. The
curve~${ E_{-17, z, w} }$ is completely defined
over~${ K = \Q(\sqrt{-17}, \sqrt{2}) }$. We compute that the trivial
character could be a splitting character for~${ E_{-17, z, w} }$,
hence the square root of the degree map could be a splitting
map~${ \beta }$. We compute that~${ c_{E_{-17, z, w}} }$
and~${ c_\beta }$ are given by
\begin{equation*}
  \begin{tabular*}{0.44\linewidth}{c|cccc}
    ${ c_{E_{-17, z, w}} }$ & ${ 1 }$ & ${ \sigma_2 }$ & ${ \sigma_{17} }$ & ${ \sigma_2 \sigma_{17} }$ \\
    \hline
    ${ 1 }$ & ${ 1 }$ & ${ 1 }$ & ${ 1 }$ & ${ 1 }$ \\
    ${ \sigma_2 }$ & ${ 1 }$ & ${ 2 }$ & ${ 1 }$ & ${ 2 }$ \\
    ${ \sigma_{17} }$ & ${ 1 }$ & ${ -1 }$ & ${ 1 }$ & ${ -1 }$ \\
    ${ \sigma_2 \sigma_{17} }$ & ${ 1 }$ & ${ -2 }$ & ${ 1 }$ & ${ -2 }$
  \end{tabular*}
  \quad \text{and} \quad
  \begin{tabular*}{0.4\linewidth}{c|cccc}
    ${ c_\beta }$ & ${ 1 }$ & ${ \sigma_2 }$ & ${ \sigma_{17} }$ & ${ \sigma_2 \sigma_{17} }$ \\
    \hline
    ${ 1 }$ & ${ 1 }$ & ${ 1 }$ & ${ 1 }$ & ${ 1 }$ \\
    ${ \sigma_2 }$ & ${ 1 }$ & ${ 2 }$ & ${ 1 }$ & ${ 2 }$ \\
    ${ \sigma_{17} }$ & ${ 1 }$ & ${ 1 }$ & ${ 1 }$ & ${ 1 }$ \\
    ${ \sigma_2 \sigma_{17} }$ & ${ 1 }$ & ${ 2 }$ & ${ 1 }$ & ${ 2 }$
  \end{tabular*}
\end{equation*}
where~${ \sigma_2 }$ and~${ \sigma_{17} }$ are generators
of~${ G_{\Q(\sqrt{2})}^K }$ and~${ G_{\Q(\sqrt{-17})}^K }$
respectively. Note that the
difference~${ c_{E_{-17, z, w}} c_\beta^{-1} }$ can not be the
coboundary of a map to~${ \Q^* }$ as it is non-symmetric, hence a
change of splitting map does not suffice. Instead, we will have to find
a twist~${ \gamma }$ for which we need a
map~${ \alpha: G_\Q^K \to K^* }$ with
coboundary~${ c_{E_{-17, z, w}} c_\beta^{-1} }$.

Using SageMath~\cite{sagemath}, we can determine
that~${ \mathcal{O}_K^* = \langle -1, \sqrt{2} - 1 \rangle }$, so
if~${ \alpha }$ would be a map with codomain~${ \mathcal{O}_K^* }$, we
can write
\begin{equation*}
  \alpha(\sigma) = (-1)^{x(\sigma)} (\sqrt{2} - 1)^{y(\sigma)}
\end{equation*}
for some~${ x, y : \G_\Q^K \to \Z }$. Now we find that
\begin{equation*}
  \begin{array}{rl}
    1 & = c_{E_{-17, z, w}} c_\beta^{-1}( \sigma_2, \sigma_{17})
    = \alpha(\sigma_2) \, \prescript{\sigma_2}{}{\alpha(\sigma_{17})} \, \alpha(\sigma_2 \sigma_{17})^{-1} \\
      & = (-1)^{x(\sigma_2) + x(\sigma_{17}) - x(\sigma_2 \sigma_{17})} \, (\sqrt{2} - 1)^{y(\sigma_2) - y(\sigma_2 \sigma_{17})} \, \prescript{\sigma_2}{}{(\sqrt{2} - 1)}^{y(\sigma_{17})} \\
    & \text{ and } \\
    -1 & = c_{E_{-17, z, w}} c_\beta^{-1}( \sigma_2, \sigma_{17})
    = \alpha(\sigma_{17}) \, \prescript{\sigma_{17}}{}{\alpha(\sigma_{2})} \, \alpha(\sigma_2 \sigma_{17})^{-1} \\
    & = (-1)^{x(\sigma_2) + x(\sigma_{17}) - x(\sigma_2 \sigma_{17})} \, (\sqrt{2} - 1)^{y(\sigma_{17}) - y(\sigma_2 \sigma_{17})} \, \prescript{\sigma_{17}}{}{(\sqrt{2} - 1)}^{y(\sigma_{2})}.
  \end{array}
\end{equation*}
Since~${ \prescript{\sigma_2}{}{(\sqrt{2} - 1)} = \sqrt{2} - 1 }$
and~${ (\sqrt{2} - 1) \, \prescript{\sigma_{17}}{}{(\sqrt{2} - 1)} =
  -1 }$, this gives us the linear system
\begin{equation*}
  \left\{
    \arraycolsep=1.4pt
    \begin{array}{rrrrrrl}
      x(\sigma_2) & + x(\sigma_{17}) & - x(\sigma_2 \sigma_{17}) &&&& \equiv 0 \text{ (mod } 2) \\
      &&& y(\sigma_2) & + y(\sigma_{17}) & - y(\sigma_2 \sigma_{17}) & = 0 \\
      x(\sigma_2) & + x(\sigma_{17}) & - x(\sigma_2 \sigma_{17}) & + y(\sigma_{2}) &&& \equiv 1 \text{ (mod } 2) \\
      &&& - y(\sigma_{2}) & + y(\sigma_{17}) & - y(\sigma_2 \sigma_{17}) & = 0,
    \end{array}
  \right.
\end{equation*}
which is clearly inconsistent. Therefore, there can be no
map~${ \alpha : G_\Q^K \to \mathcal{O}_K^* }$ with
coboundary~${ c_{E_{-17, z, w}} c_\beta^{-1} }$.

Now look at the map~${ \alpha : G_\Q^K \to K^* }$ given by
\begin{equation*}
  \alpha(\sigma) = \left\{
    \begin{array}{cl}
      1 & \text{if } \sigma = 1 \\
      -1 & \text{if } \sigma = \sigma_2 \\
      \frac{1 - 3 \sqrt{2}}{\sqrt{-17}} & \text{otherwise.}
    \end{array}
  \right.
\end{equation*}
A quick computation shows that the coboundary of this map
is~${ c_{E_{-17, z, w}} c_\beta^{-1} }$, and
that~${ \gamma = 1 + 3 \sqrt{2} }$
satisfies~${ \prescript{\sigma}{}{\gamma} = \alpha(\sigma)^2 \gamma }$
for all~${ \sigma \in G_\Q^K }$.
\begin{remark}
  The codomain of this~${ \alpha }$ is the ring of~${ S
  }$-units~${ \mathcal{O}_{K, S}^* }$, where~${ S }$ is the set of all
  primes above~${ 17 }$.
  As shown in~\cite[Proposition 2.5.5]{LangenThesis},
  one can always compute a finite set of primes~${ S }$ such that
  there exists an~${ \alpha }$ with image in~${ \mathcal{O}_{K, S}^* }$.
  This example shows that in some cases a non-empty~${ S }$ is actually
  necessary.
\end{remark}

\begin{remark}
  Note that 
  \begin{equation*}
    1 - \sqrt{-17} = (1 + 3 \sqrt{2}) \left(\frac{1}{\sqrt{2}} + \frac{3}{\sqrt{-17}} - \frac{1}{\sqrt{2} \sqrt{-17}} \right)^2,
  \end{equation*}
  so we could also take~${ \gamma = 1 - \sqrt{-17} }$. The advantage
  of this is that~${ E_{-17, z, w}^\gamma }$ remains defined
  over~${ \Q( \sqrt{-17} ) }$ and in fact it is even completely
  defined over that field. The problem is that~${ 1 - \sqrt{-17} }$ is
  divisible by a prime above~${ 3 }$ in~${ \Q(\sqrt{-17}) }$.  By
  Proposition~\ref{thm:levelLowering} a~${ 3 }$ will appear in the
  level after level lowering, which most likely increases the
  dimension and thus computation time of the space of newforms.
\end{remark}

The curve~${ E_{-17, z, w}^{\gamma} }$ remains completely defined
over~${ K }$. The conductor of~${ E_{-17, z, w}^{\gamma} }$
over~${ K }$ is
\begin{equation*}
  \left\{
    \begin{array}{cl}
      \mathfrak{p}_2^{16} \, (17) \left(\Rad_{2\cdot17}(w^2 + 17 \, z^4)\right) & \text{if } 2 \mid w \\
      \mathfrak{p}_2^{6} \, (17) \left(\Rad_{2\cdot17}(w^2 + 17 \, z^4)\right) & \text{if } 2 \nmid w.
    \end{array}
  \right.
\end{equation*}
Here~${ \mathfrak{p}_2 }$ is the unique prime of~${ K }$
above~${ 2 }$. The field generated by the values of~${ \beta }$
is~${ \Q(\sqrt{2}) }$. Since this is quadratic, we know
that~${ \residue_\Q^{K} E_{-17, z, w}^{\gamma} }$ is a product of
two~${ \Q }$-simple abelian variety of~${ \GL_2 }$-type. The levels of
the associated newforms can be computed to be
\begin{equation*}
  \left\{
    \begin{array}{ll}
      (2^8 \cdot 17^2 \cdot \Rad_{2\cdot17}(w^2 + 17 \, z^4), 2^8 \cdot 17^2 \cdot \Rad_{2\cdot17}(w^2 + 17 \, z^4)) & \text{if } 2 \mid w \\
      (2^5 \cdot 17^2 \cdot \Rad_{2\cdot17}(w^2 + 17 \, z^4), 2^6 \cdot 17^2 \cdot \Rad_{2\cdot17}(w^2 + 17 \, z^4)) & \text{if } 2 \nmid w.
    \end{array}
  \right.
\end{equation*}
So after level lowering, the lowest possible levels
are~${ 2^8 \cdot 17^2 }$ if~${ 2 \mid w }$ and~${ 2^5 \cdot 17^2 }$
if~${ 2 \nmid w }$. Note that the character of these newforms is
trivial as the splitting character for all splitting maps is trivial.

Note that we can apply level lowering when the corresponding
mod~${ l }$ Galois representation is irreducible. This is the case
for~${ l = 3, 5, 7, 13 }$ by
Corollary~\ref{thm:freycurveirreducible}. For~${ l = 11 }$
or~${ l > 13 }$ we need that the corresponding
point~${ P \in E_D(\Q) \setminus \{ \mathcal{O}, T \} }$ has~${ B_P }$
divisible by a prime number~${ p > 3 }$. The points for which this is
not the case can be computed to be
\begin{equation*}
  \left\{
    \begin{array}{c}
      \pm P, \pm Q, \pm 2 P, \pm 2 Q, \pm P + T, \pm Q + T, \pm 2 Q + T, \pm ( P + Q ), \\
      \pm (P - Q), \pm (2 P - 2 Q), \pm (P - Q) + T, \pm (P - 2 Q) + T
    \end{array}
  \right\}.
\end{equation*}
The only perfect powers among the~${ B }$ of these points
are~${ B_{\pm 2 P} = B_{\pm 2 Q} = 2^2 }$
and~${ B_{\pm 2 Q + T} = 3^2 }$. We may thus assume that for odd prime
numbers~${ l }$ we have irreducibility.

Computing the newforms of these levels takes a few hours with
Magma. Using the framework, we find that by comparing traces at the
primes
\begin{equation*}
  \{ 3, 5, 7, 11, 13, 19, 29, 31, 37, 41, 43, 47, 59, 67, 73, 97, 113 \},
\end{equation*}
we can eliminate all but ${ 10 }$ newforms for all primes~${ l > 31 }$.
\begin{remark}
  We have performed the elimination for all prime numbers below 200
  besides~${ 2 }$ and~${ 17 }$, but the list given here includes the
  primes at which elimination actually happened.
\end{remark}

The~${ 10 }$ newforms that remain, all seem to have complex
multiplication and coefficient field~${ \Q(\sqrt{2}) }$. They most
likely correspond to values of~${ (z, w) }$ corresponding to (pseudo)
solutions of the problem. For example, the point~${ T }$ corresponds
to the curve~${ E_{-17, 0, 1}^\gamma }$, which in turn corresponds to
one of the newforms of level~${ 2^5 \cdot 17^2 }$.

To avoid these (pseudo-)solutions, we look at non-zero
multiples~${ P_m }$ of the point~${ P_1 = P + Q + T }$.
Since~${ B_1 = 7 }$, \eqref{eqn:edsp} tells us that all~${ B_m }$ are
divisible by~${ 7 }$. Using this additional information when comparing
traces of Frobenius at~${ 7 }$, allows us to eliminate all newforms
for~${ l > 17 }$. This leads to the following asymptotic result.

\begin{theorem}
  \label{thm:Dm17}
  The sequence~${ B_m }$ contains no~${ l }$-th powers with~${ l }$ a
  prime number
  when~${ P_1 = \left(-\frac{153}{49}, \pm\frac{1632}{343}\right) \in
    E_{-17}(\Q) }$ and~${ l > 17 }$.
\end{theorem}

\begin{remark}
  Besides the methods described here, we also tried to eliminate the
  newforms~${ l = 13 }$ and~${ l =17 }$ by performing Kraus'
  method. This did not allow us to eliminate all newforms for these
  primes. 
\end{remark}

\subsection{Further examples}\label{FurtherExamples}
Using the framework, one can easily compute further asymptotic results
as in Theorem~\ref{thm:D125} and Theorem~\ref{thm:Dm17}. We
will not write this out in detail as the approach is very
similar. Table~\ref{tab:EDSresults} contains the results of the computations
that were performed.

\begin{table}[h!]
  \centering
  \renewcommand{\arraystretch}{1.3}
  \begin{tabular}{c|c|c|c}
    case & $D$ & $P_1 \in E_D(\Q)$ & $l$ \\
    \hline
    (i) & $125$ & $\left( \frac{121}{4}, \pm \frac{1\,419}{8} \right)$ & $ l > 2 $ \\
    (ii) & $125$ & $\left( \frac{500}{121}, \pm \frac{32\,250}{1\,331} \right) = \left( \frac{121}{4}, \mp \frac{1\,419}{8} \right) + (0,0)$ & $ l > 5 $, $ l \ne 11 $ \\
    (iii) & $-17$ & $\left(-\frac{153}{49}, \pm\frac{1\,632}{343}\right)$ & $ l > 17 $ \\
    (iv) & $3$ & $ \left(\frac{1}{4}, \pm \frac{7}{8} \right)=2(3,\pm 6) $ & $ l > 2 $ \\
    (v) & $3$ & $ \left(\frac{27}{121}, \pm\frac{1\,098}{1\,331} \right)=3(3,\mp 6) $ & $ l > 17 $ \\
    (vi) & $-2$ & $ \left(\frac{9}{4}, \pm\frac{21}{8}\right)=2(-1,\mp 1) $ & $ l > 2 $ \\
    (vii) & $-2$ & $ \left(-\frac{1}{169}, \pm \frac{239}{2197}\right)=3(-1,\pm 1) $ & $ l > 2 $ \\
    (viii) & $-2$ & $ \left(\frac{4\,651\,250}{1\,803\,649}, \pm \frac{8\,388\,283\,850}{2\,422\,300\,607}\right)=5(-1,\mp 1)+(0,0) $ & $ l > 5 $, $ l \ne 79 $ \\
    (ix) & $-2$ & $ \left(-\frac{8}{9}, \pm \frac{28}{27}\right)=2(-1,\mp 1)+(0,0) $ & $ l > 3 $ \\
  \end{tabular}
  \caption{Elliptic divisibility sequences with no $l$-th powers}
  \label{tab:EDSresults}
\end{table}

\subsection{Small exponent values}

While our main focus in this article is bounding prime exponents, a natural next step would be to determine all perfect powers in the elliptic divisibility sequences we consider. It follows immediately from~\cite[Theorem 1.1]{EverestReynoldsStevens07} (dealing with any nonsingular Weierstrass equation over $\Z$) that for every integer $l>1$, our main equation
\begin{equation}~\label{eqn:main}
B_m=v^l, \quad m,v \in \Z_{>0}
\end{equation}
has only finitely many solutions. In particular, it suffices to restrict to prime exponents $l$.
We note that this finiteness result is not effective, as it appeals to Faltings' finiteness theorem for rational points on curves of genus greater than one, amongst other things.
Using the fact that our elliptic curves $E_D$ are of a rather special form~\eqref{eqn j=1728 curve}, we will now discuss an independent, rather direct, reduction of solving~\eqref{eqn:main} for a fixed integer $l>1$ to finding $\Q$-rational points on finitely many hyperelliptic curves over $\Q$ of genus $2l-1>1$ (or hyperelliptic quotients thereof of smaller genus). In general this approach would lead again, by Faltings' theorem, to an ineffective finiteness result for fixed exponents in~\eqref{eqn:main}. But in favourable cases it could fall within the scope of effective methods for determining rational points on (hyperelliptic) curves, though we shall not investigate this much in this article.

Fix some integers $l>1$. For the construction, we start by recalling that~\eqref{eqn basic 3-term}, with $\hat{B}=:z$ and $B=B_m=v^l$, leads to a generalized Fermat equation of signature $(2,4,4l)$, namely
\begin{equation}\label{eqn:GFE 4,4,4l}
w^2=az^4+\hat{a}v^{4l},
\end{equation}
to be solved in pairwise coprime (positive) integers $w,z,v$.
We remark that the sum of the reciprocals of the exponents satisfy $1/2+1/4+1/(4l)<1$, so by~\cite[Theorem 2]{DarmonGranville95} there are only finitely many solutions, where once again the proof reduces the finiteness to Faltings' theorem.
Now instead of considering~\eqref{eqn:GFE 4,4,4l} of signature $(2,4,l)$, which lead to our Frey $\Q$-curve construction, we will consider it of signature $(2,2,4l)$, i.e., setting $u:=z^2$, we get
\begin{equation}\label{eqn:GFE 2,2,4l}
w^2=au^2+\hat{a}v^{4l}
\end{equation}
to be solved in pairwise coprime (positive) integers $w,u,v$.
Since $1/2+1/2+1/(4l)>1$, it is a \emph{spherical} generalized Fermat equation, which yields that the solutions are given by finitely many parametrizations.
More precisely, there exist finitely many, say $r$, triples of separable binary forms $F_i, G_i, H_i \in \Z[r,s]$ ($i=1,2,\ldots,r$) of degrees $4l, 4l, 2$ respectively, satisfying $F_i^2=aG_i^2+\hat{a} H_i^{4l}$, and such that for any pairwise coprime integer solution $(w,u,v)$ to~\eqref{eqn:GFE 2,2,4l} there exists an index $i$ and coprime integers $r,s$ such that $(F_i(r,s),G_i(r,s),H_i(r,s))=(w,u,v)$. This means that a pariwise coprime solution $(w,z,u)$ to the original equation~\eqref{eqn:GFE 4,4,4l} satisfies $z^2=G_i(r,s)$ for some index $i$.
Each of these equations defines a hyperelliptic curve 
\begin{equation}\label{eqn:hyperelliptic}
C_i: z^2=G_i(r,s)
\end{equation}
over $\Q$ in weighted projective space (of weights $1,1,2l$ for $r,s,z$ respectively) of genus $2l-1$ (since the $G_i$ are separable). The $C_i$ have some interesting quotients that might be helpful in determining the rational points $C_i(\Q)$.

\begin{example}
We look at the elliptic divisibility sequence from Subsection~\ref{subsec:example 125}, where $D=125$ and we choose $P=(121/4,\pm 1419/8)$. We note that we have $a=1$, and hence $\hat{a}=125$, in~\eqref{eqn:GFE 4,4,4l} for solutions to~\eqref{eqn:main}. We can perform an elementary descent over $\Z$ to obtain sufficiently many hyperelliptic curves $C_i$ as in~\eqref{eqn:hyperelliptic} by rewriting~\eqref{eqn:GFE 4,4,4l} as $(w+z^2)(w-z^2)=5^3v^{4l}$. Recall that we have $2|B_m$. Hence $2|v$, which leads to $\gcd(w+z^2,w-z^2)=2$. Without loss of generality we can and will assume $w>0$, which now leads to
\[w+z^2=2 c_2 c_5 \alpha^{4l}, \qquad w-z^2=2 c_2' c_5' \beta^{4l}\]
for some $\alpha, \beta \in \Z$ and $\{c_2,c_2'\}=\{1,2^{4l-2}\}, \{c_5,c_5'\}=\{1,5^3\}$.
Subtracting the second equation from the first and dividing by $2$ yields equations for our hyperelliptic curves:
\begin{equation}
C_i: z^2=c_2 c_5 \alpha^{4l}-c_2' c_5' \beta^{4l}
\end{equation}
for $i=1,2,3,4$, say by choosing $(c_2,c_5)=(1,1), (1,5^3), (2^{4l-2},1), (2^{4l-2},5^3)$ respectively (which then also fixes the corresponding $(c_2',c_5')$).

As indicated by Theorem~\ref{thm:D125}, the only (prime) exponent left to deal with is $l=2$, so let us fix this value for the rest of this example. The hyperelliptic curves $C_i$ are of genus $3$ and one easily obtains that $C_2(\Q)=\emptyset$ by checking locally that $C_2(\Q_2)=\emptyset$. On all of the other 3 curves one can easily spot some $\Q$-rational point. All curves have genus $2$ quotients, of which the jacobians all turn out to have rank $2$, so that Chabauty-Coleman does not immediately apply to find all rational points. The equations invite some further descent, but we will not pursue determining $C_i(\Q)$ for $i=1,3,4$ further here.

As a final note, equations for the $l=2$ case can also be obtained by considering~\eqref{eqn:GFE 4,4,4l} of signature $(2,4,2)$, whose solutions can again be parametrized. This leads to finitely many binary quartic forms $E_i \in \Z[s,t]$ such that solutions to the original equation are given by $v^4=E_i(s,t)$, which define projective plane quartic curves.

\end{example}

\subsection{Alternative approaches for some examples}
\label{sec:EDSalternatives}
Consider any elliptic curve $E/\Q$ (in Weierstrass form with integral coefficients) and non-torsion point $P=(x,y) \in E(\Q)$ with the denominator of $x$ divisible by $p \in \{2,3\}$.
Reynolds associates in~\cite{Reynolds12} a Frey elliptic curve over $\Q$ (depending on $p$) to our Diophantine problem of interest~\eqref{eqn:main}.
Together with Silverman's famous result on the existence of primitive divisors~\cite{Silverman88} it is then shown that there exists an effective bound $l_0$ such that for all primes $l \geq l_0$ there are no solutions to~\eqref{eqn:main}; see Theorem 1.2 in \emph{loc. cit.} (where again it is also noted that consequently $(B_m)$ contains only finitely many perfect powers).
The cases (i), (iv), (vi), and (ix) from Table~\ref{tab:EDSresults} fall in this category. The Frey curve mentioned above is actually associated to the elliptic divisibility subsequence for $pP$. As such, for cases (v) and (vii) from Table~\ref{tab:EDSresults}, there is also an alternative approach using Frey curves over $\Q$. This leaves cases (ii), (iii), and (viii) for which there does not seem to be alternative approaches available in the literature.

In the case that~${ a }$ is a positive non-square, the Frey
curve~${ E_{a, z, w} = E_{a, z, w}^1 }$ is defined over a totally real
field~${ \Q(\sqrt{a}) }$. As an alternative to the~${ \Q }$-curve
approach, we could therefore use modularity of elliptic curves over
totally real quadratic fields \cite{FreitasHungSiksek15} and perform
the modular method with the associated Hilbert modular
forms, as already mentioned (for the subcase $D>0$) in the introduction.
Furthermore, one could try to combine the~${ \Q }$-curve
approach with this approach as a multi-Frey method.

We have tried to apply this Hilbert modular approach to some
examples. In the case~${ a = D = 3 }$, this did not lead to any
additional information. For larger~${ D }$, the levels of the Hilbert
modular forms became too large to feasibly compute the corresponding
newforms.


\bibliographystyle{amsplain}
\bibliography{BoundingPPinEDS}

\end{document}